\documentclass[11pt, reqno]{amsart}

\usepackage{amssymb}
\usepackage[T1]{fontenc}

\evensidemargin\oddsidemargin

\begin{document}
\baselineskip=18pt
\setcounter{page}{1}
    
\newtheorem{Conj}{Conjecture}
\newtheorem{TheoA}{Theorem A\!\!}
\newtheorem{TheoB}{Theorem B\!\!}
\newtheorem{Lemm}{Lemma}
\newtheorem{Rem}{Remark}
\newtheorem{Coro}{Corollary}
\newtheorem{Propo}{Proposition}

\renewcommand{\theTheoA}{}
\renewcommand{\theTheoB}{}

\def\a{\alpha}
\def\b{\beta}
\def\B{{\bf B}} 
\def\C{{\bf C}} 
\def\cG{{\mathcal{G}}} 
\def\cH{{\mathcal{H}}} 
\def\cI{{\mathcal{I}}} 
\def\cS{{\mathcal{S}}}
\def\UU{{\mathcal{U}}}
\def\ca{c_{\a}}
\def\ka{\kappa_{\a}}
\def\coa{c_{\a, 0}}
\def\cua{c_{\a, u}}
\def\cL{{\mathcal{L}}} 
\def\cM{{\mathcal{M}}} 
\def\Ea{E_\a}
\def\eps{{\varepsilon}} 
\def\esp{{\mathbb{E}}} 
\def\Ga{{\Gamma}} 
\def\G{{\bf \Gamma}} 
\def\GG{{\bf G}}
\def\HH{{\bf H}}
\def\ii{{\rm i}}
\def\e{{\rm e}}
\def\L{{\bf L}}
\def\lbd{\lambda}
\def\lacc{\left\{}
\def\lcr{\left[}
\def\lpa{\left(}
\def\lva{\left|}
\def\M{{\bf M}}
\def\NN{{\mathbb{N}}} 
\def\pb{{\mathbb{P}}}
\def\pab{{\varphi_{a,b}}} 
\def\tpab{{{\widetilde \varphi}_{a,b}}} 
\def\rl{{\mathbb{R}}}
\def\racc{\right\}}
\def\rpa{\right)}
\def\rcr{\right]}
\def\rva{\right|}
\def\prost{{\succ_{\! st}}}
\def\W{{\bf W}}
\def\X{{\bf X}}
\def\XX{{\mathcal X}}
\def\Y{{\bf Y}}
\def\U{{\bf U}}
\def\V{{\bf V}_\a}
\def\Un{{\bf 1}}
\def\Z{{\bf Z}}
\def\A{{\bf A}}
\def\AA{{\mathcal A}}
\def\hAA{{\hat \AA}}
\def\hL{{\hat L}}
\def\hT{{\hat T}}

\def\claw{\stackrel{d}{\longrightarrow}}
\def\elaw{\stackrel{d}{=}}
\def\qed{\hfill$\square$}

\newcommand*\pFqskip{8mu}
\catcode`,\active
\newcommand*\pFq{\begingroup
        \catcode`\,\active
        \def ,{\mskip\pFqskip\relax}%
        \dopFq
}
\catcode`\,12
\def\dopFq#1#2#3#4#5{%
        {}_{#1}F_{#2}\biggl[\genfrac..{0pt}{}{#3}{#4};#5\biggr]%
        \endgroup
}

\title{Convolution of beta prime distribution}

\author[Rui. A. C. Ferreira]{Rui A. C. Ferreira}

\address{Grupo F\'isica-Matem\'atica, Faculdade de Ci\^encias, Universidade de Lisboa,
Av. Prof. Gama Pinto 2, 1649-003 Lisboa, Portugal. {\em Email}: {\tt raferreira@fc.ul.pt}}

\author[Thomas Simon]{Thomas Simon}

\address{Universit\'e de Lille, CNRS, UMR 8524 - Laboratoire Paul Painlevé, 59000 Lille, France. 
{\em Email}: {\tt thomas.simon@univ-lille.fr}}



\begin{abstract} 
We establish some identities in law for the convolution of a beta prime distribution with itself, involving the square root of beta distributions. The proof of these identities relies on transformations on generalized hypergeometric series obtained via Appell series of the first kind and Thomae's relationships for ${}_3F_2(1)$. Using a self-decomposability argument, the identities are applied to derive complete monotonicity properties for quotients of confluent hypergeometric functions having a doubling character. By means of probability, we also obtain a simple proof of Tur\'an's inequality for the parabolic cylinder function and the confluent hypergeometric function of the second kind. The case of Mill's ratio is discussed in detail.
\end{abstract}

\maketitle

\section{Introduction}

The beta prime random variable $\B'_{a,b}$ with parameters $a,b > 0$ is a positive random variable with density
$$\frac{\Ga(a+b)}{\Ga(a)\Ga (b)} \,\frac{x^{a -1}}{(1+x)^{a+b}}$$
on $(0,\infty).$ This distribution is well-known in statistics as a generalization of the classical Fisher-Snedecor distribution: if $d_1, d_2$ are positive integers and $\chi^{}_2(d_1), \chi^{}_2(d_2)$ are two independent chi-squared random variables with respective parameters $d_1,d_2$, then the ratio
$$\frac{\chi^{}_2(d_1)}{\chi^{}_2(d_2)}$$
which is used in testing the equality of variances, is distributed as $\B'_{\frac{d_1}{2},\frac{d_2}{2}}.$ In general, one has the identities
\begin{equation}
\label{base}
\B'_{a,b}\; \elaw\; \frac{1}{\B_{b,a}} \, -\, 1\; \elaw\; \frac{\G_a}{\G_b} 
\end{equation}
where $\B_{p,q}$ and $\G_t$ denote for every $p,q,t >0$ the standard beta and gamma random variables with respective densities
$$\frac{\Ga(p+q)}{\Ga(p)\Ga (q)} \, x^{p -1}(1-x)^{q-1}\,\Un_{(0,1)}(x)\qquad\mbox{and}\qquad \frac{1}{\Ga(t)} \, x^{t -1} e^{-x}\,\Un_{(0,\infty)}(x)$$
and, here and throughout, all random variables appearing on one side of an identity in law are assumed to be independent unless otherwise explicitly stated. In the statistical literature, the beta prime distribution is sometimes called the beta distribution of the second kind, the generalized $F-$distribution or the Pearson type VI distribution, and we refer to Chapter 27 in \cite{BJK} for an account.

The Mellin transform of the beta prime distribution is given as
\begin{equation}
\label{Mell}
\esp[(\B'_{a,b})^s] \; =\; \frac{\Ga(a+s)\Ga(b-s)}{\Ga(a)\Ga(b)}
\end{equation}
for every $s\in (-a,b),$ whereas its Laplace transform reads
\begin{equation}
\label{LT}
\esp[e^{-z \B'_{a,b}}] \; =\; \frac{\Ga(a+b)}{\Ga(a)\Ga (b)} \int_0^\infty e^{-z t} t^{a -1} (1+t)^{-a-b}\, dt \; =\; \frac{\Ga(a+b)}{\Ga (b)}\, \Psi(a,1-b,z)
\end{equation}
for every $z \ge 0$, where $\Psi$ is the confluent hypergeometric function of the second kind - see Formula 6.5.(2) in \cite{EMOT}. This connection with the $\Psi$ function was used by Ismail and Kelker in \cite{IK} to show, by Stieltjes inversion, the non-obvious fact that the beta prime distribution is self-decomposable (SD), in other words that for every $a,b >0$ and $\lbd\in (0,1),$ there is a decomposition
$$\B'_{a,b}\;\elaw\; \lbd\,\B'_{a,b}\; +\; \X_{a,b,\lbd}$$
for some positive, infinitely divisible (ID) random variable $\X_{a,b,\lbd}.$ This property is characterized by the following L\'evy-Khintchine representation of $\Psi(a,c,z)$ which is a consequence of (1.4) in \cite{IK} by integration:
$$\Psi(a,c,z)\; =\; \frac{\Ga(1-c)}{\Ga(a+1-c)}\, \exp\, -\lpa\int_0^\infty (1- e^{-zt})\, t^{-1}\, k_{a,c}(t)\, dt\rpa$$
for every $a,z>0$ and $c<1,$ with 
$$k_{a,c}(t)\; =\; \int_0^\infty e^{-ty}\lpa \frac{y^{-c} e^{-y} \lva \Psi(a,c, y e^{\ii \pi})\rva^{-2}}{\Ga(a)\Ga(a+1-c)}\rpa dy$$
a decreasing function on $(0,\infty).$ In probabilistic terms, this formula means that the Lévy measure of the self-decomposable random variable $\B'_{a,b}$ is expressed in terms of the values on the cut of the same confluent hypergeometric function $\Psi(a,1-b,z).$ We refer e.g. to \cite{Bond} p.18 for a brief, enough for our purposes in the present paper, account on self-decomposability. Throughout, we will use the classical fact that this class of distributions is closed w.r.t. convolution and weak convergence.

The L\'evy-Khintchine representation means that $\Psi(a,c,z)$ is for $a >0$ and $c<1$ not only completely monotone (CM) as the renormalized Laplace transform of a beta prime random variable, but also logarithmically completely monotone (LCM), that is the logarithmic derivative $z\mapsto -(\log \Psi(a,c,z))'$ is CM. Throughout the paper, the CM property will be implicitly meant on $(0,\infty).$ The completely monotonic character of the spectral function $k_{a,c}$ itself means a stronger property for the beta prime distribution which is called the generalized gamma convolution and amounts to the stochastic integral representation
$$\B'_{a,b}\; \elaw\; \int_0^a\! \lbd_{a,b}(t)\, d\gamma_t,$$
where $\{\gamma_t,\; t\ge 0\}$ is the standard Gamma subordinator and $\{\lbd_{a,b}(t), \; t\in[0,a]\}$ some deterministic function that can be recovered from $k_{a, 1-b}.$ Notice that the upper integration bound $a = k_{a,1-b} (0)$ corresponds to the power law behaviour in $x^{a-1}$ at zero for the density of $\B'_{a,b}.$ In this respect, the function $t\mapsto a^{-1}k_{a,1-b}(t)$ is the Laplace transform of a positive random variable $\GG_{a,b}$ which is called the Thorin random variable associated to $\B'_{a,b}$. The distribution of $\GG_{a,b}$ is quite interesting and will be reviewed in Section 2. We refer to chapters 3 and 4 in \cite{Bond} and Section 1 in \cite{JRY} for accounts on generalized Gamma convolutions and Thorin random variables, including the proofs of the above facts. We also mention that around the same time when \cite{IK} was written, it was observed by Bondesson and Thorin that the density of $\B'_{a,b}$ is the prototype of a so-called hyperbolically completely monotone (HCM) density, which has stronger properties than the densities of generalized Gamma convolutions. We refer to Chapter 5 in \cite{Bond} for more information on this notion, which we will briefly use in Section 5.    

In this paper, we investigate the structure of the addition of two independent beta prime random variables. This problem, which is naturally connected to the semigroup property underlying the random variable $\B'_{a,b}$, was investigated in \cite{D} with motivations coming from Bayesian analysis and multivariate testing. In particular, a general and rather complicated formula for the distribution function is given in \cite{D} in terms of some Lauricella function in three variables, which was later simplified in \cite{E} into an Appell series in two variables. Our main purpose here is to display some simple, albeit non-trivial, identities in law having a rather surprising character involving a product, a sum and a square root. Our first result is the following.
  
\begin{TheoA} For every $a > 0,$ one has
$$\B'_{a,1/2} \, + \, \B'_{a,1/2} \; \elaw \; \B'_{2a, 1/2} \,\times \lpa 1 \, +\, \sqrt{\B_{a,1/2}}\rpa.$$
\end{TheoA}

This identity, which should be compared with the well-known identities
$$\frac{1}{\G_{1/2}}\; +\; \frac{1}{\G_{1/2}}\; \elaw\; \frac{4}{\G_{1/2}}\qquad\quad\mbox{and}\qquad\quad \G_a\; +\; \G_a\; \elaw\; \G_{2a},$$ 
amounts to the computation of the law of a certain ``Dirichlet mean'' associated to the Thorin random variable $\GG_{a,1/2}:$ with the notation of \cite{JRY}, to which we refer for details on Dirichlet means, one has
$$D_a(\GG_{a,1/2})\; \elaw\; \frac{1}{\G_{1/2}}\qquad\mbox{and}\qquad D_{2a}(\GG_{a,1/2})\; \elaw\; \frac{1+\sqrt{\B_{a,1/2}}}{\G_{1/2}}\cdot$$ 
Apart from its exactitude, our motivation for the above result comes from sharp bounds, monotonicity and complete monotonicity properties for quotients of confluent hypergeometric functions. The complete monotonicity of the ratio
$$\frac{\Psi(a+1,c+1,z)}{\Psi(a,c,z)}$$
for every $a >0$ and $c < 1,$ which is equivalent to the infinite divisibility of $\B'_{a,b}$ for every $a,b > 0,$ was established in \cite{IK} and has influenced ever since a substantial body of works devoted to such questions for quotients of various classes of special functions. We refer to \cite{BI, BPS, IL, K, LN, MS, Seg1, Seg3, TS, YZ} for a few examples of this literature, which is still ongoing. The approach which is generally followed in these papers is analytical and often relies on the second order ODE associated to the class of special functions in consideration. As in the aforementioned result of \cite{IK}, this leads to a discrete difference between the parameters of the two special functions appearing in the quotient. This discrete gap is also observed in the so-called Tur\'an inequalities, an emblematic example of which is the recently established - see \cite{K} - monotonicity on $\rl$ of the quotient of parabolic cylinder functions  
$$\frac{(D_{-\nu-1}(z))^2}{D_{-\nu}(z) D_{-\nu -2}(z)}$$
for every $\nu > 0.$ In the present paper, we take a probabilistic point of view and show properties of this kind with a difference in the parameters which is continous. Typical results that we obtain are the complete monotonocity of   
$$\frac{\Psi(a,c,z)}{\Psi(a,c',z)}$$
for every $a >0$ and $c' < c < 1$ which is proved by self-decomposability - see Corollary \ref{B1}, and the monotonicity on $\rl$ of
$$\frac{(D_{-\nu-c}(z))^2}{D_{-\nu}(z) D_{-\nu -2c}(z)}$$
for every $\nu, c > 0,$ which we derive through the connection with generalized half-Gaussian distributions - see Proposition \ref{Hermit}. The self-decomposability argument also implies as a by-product of Theorem A the complete monotonicity of the quotient of Hermite function
$$\frac{(H_{-\nu}(\sqrt{z}))^2}{H_{-2\nu}(\sqrt{z})}$$
for every $\nu > 0,$ whose ``doubling character'' is less usual than the previous results. For $\nu\in (0,2]$, a stochastic comparison of the distribution of the Dirichlet means $D_a(\GG_{a,1/2})$ and $D_{2a}(\GG_{a,1/2}),$ and a general formula due to \cite{JRY} imply that the above quotient is actually LCM, in other words that its associated random variable is compound Poisson - see Corollary \ref{Para2}. In the last part of this paper, these results allow us to derive a series of properties of the classical Mill's ratio, including a new proof of Sampford's inequality.\\

Our second exact result for the convolution of a beta prime distribution is the following identity.

\begin{TheoB} For every $b\in (0,1/2),$ one has
$$\B'_{a,b} \, + \, \B'_{a,b} \; \elaw \; \B'_{2a, b} \,\times \lpa 1 \, +\, \sqrt{\frac{\B_{a,1/2}}{\B_{b,1/2-b}}}\rpa$$
if $a = 1 - b$ or $a =1/2.$
\end{TheoB}

In the case $a = 1/2,$ this result recovers Theorem A in letting $b\to 1/2$ and $\B_{b,1/2-b} \claw 1.$ The appearance of a third random variable on the right-hand side of the identity makes Theorem B significantly harder to prove than Theorem A. Whereas the latter relies on rather classical considerations on the hypergeometric series ${}_2F_1(z),$ the analysis underlying the proof of the former makes an extensive use of the generalized series ${}_3F_2(z),$ where the corresponding formul\ae \, are more sparse. Our argument is based on certain transformations of the Mellin transform of the left-hand side, where Thomae's relationships on ${}_3F_2(1)$ play their role, and the comparison of certain densities obtained via some non-trivial reduction formul\ae\, on Appell series of the first kind. A key-argument is also the transformation of a certain ${}_3F_2(1)$ into some ${}_3F_2(-1),$ which is necessary to obtain a random variable of the type $\B'_{2a,b}\times(1+\X)$ on the right-hand side. Unfortunately, this transformation is possible only in certain specific situations and this explains our limitation on the range of the parameters. As a further, negative, result, we prove by stochastic ordering that for every $a > 0,$ an identity in law with a random variable of the type $\B'_{2a,b}\times (1+\X)$ on the right-hand side cannot hold if $b > 1$ - see Proposition \ref{Stoo}. On the other hand, we do believe that the statement of Theorem B is valid for every $a > 0$ and $b\in (0,1/2)$ and it is our hope that this probabilistic problem, which is equivalent to a non-conventional integral representation for a family of ${}_3F_2(1)$'s, can stimulate readers having a deeper grip on generalized hypergeometric series. 

\section{Preliminaries}

In this section, we exploit an elementary identity in law stated as Proposition \ref{B0}, in order to derive a CM property for quotients of $\Psi$ functions by a self-decomposability argument. The same argument will be used in the next section to deduce some further CM properties from the more involved Theorems A and B. In the case $a\le 1,$ we also obtain a closed and non-trivial formula for the cumulative Thorin measure of $\B'_{a,b}$ which is used to derive the stronger LCM property. In passing, we establish a connection between this Thorin measure and the so-called Askey-Wimp-Kerov distribution. 

\begin{Propo} 
\label{B0}
For every $b'> b > 0$ and $a > 0,$ one has
$$\B'_{a,b} \; \elaw\; \B'_{a,b'} \,\times \lpa 1\, +\, \B'_{b'\!-b,b}\rpa$$
\end{Propo}

\proof
The factorization 
\begin{equation}
\label{base1}
\G_b\; \elaw\; \G_{b'}\, \times \,\B_{b, b'\!-b}
\end{equation}
is well-known and follows from the computation of the Mellin transform on both sides. Plugging this into (\ref{base}) implies  
$$\B'_{a,b}\; \elaw\; \B'_{a,b'}\,\times \,\B_{b,b'\!-b}^{-1} \; \elaw\; \B'_{a,b'} \, \times \lpa 1\, +\, \B'_{b'\!-b,b}\rpa$$
as required.
\endproof

As a consequence of this identity in law, we obtain our first complete monotonicity property. 
 
\begin{Coro}
\label{B1}
For every $a > 0$ and $c' < c < 1,$ the function
$$z\; \mapsto\; \frac{\Psi(a, c, z)}{\Psi(a, c', z)}$$  
is {\em CM}.
\end{Coro}

\proof
Set $b = 1- c, b' = 1- c'$ with $b' > b > 0.$ As discussed in the introduction, the random variable $\B'_{a,b'}$ is SD, and this implies that conditionally on $\B'_{b'\!-b,b} = x \ge 0,$ the law of the independent product $\B'_{a,b'} \, \times (1\, +\, \B'_{b'\!-b,b}) = (1 + x)\, \B'_{a,b'}$ is that of the independent sum
$$\B'_{a,b'} \, +\, \X^{a,b'}_{x}$$
for some non-negative random variable $\X^{a,b'}_x\!\!.$ Integrating shows, by Proposition \ref{B0}, the additive factorization
$$\B'_{a,b} \; \elaw\; \B'_{a,b'}\; +\; \Y_{a,b,b'}$$
with the notation $\X^{a,b'}_{\B'_{b'\!-b,b}} \!\!\! = \Y_{a,b,b'}\ge 0.$ By (\ref{LT}), we deduce the representation
$$\frac{\Psi(a, c, z)}{\Psi(a, c', z)}\; =\; \frac{\Ga(a+b)\Ga(b')}{\Ga(a+b')\Ga (b)}\;\esp\lcr e^{-z \Y_{a,b,b'}}\rcr,$$  
which concludes the proof by Bernstein's theorem.

\endproof

\begin{Rem} 
\label{Kumma}
{\em For every $a > 0$ and $1 < c < c' < 1+a,$ Kummer's transformation - see Formula 6.5.(6) in \cite{EMOT} - implies that
$$z\;\mapsto\; z^{c-c'} \, \frac{\Psi(a + c-c', c, z)}{\Psi(a, c', z)}\; =\; \frac{\Psi(1+a -c', 2-c, z)}{\Psi(1+a -c', 2-c', z)}$$
is also CM. The complete monotonicity of
$$z\;\mapsto\;\frac{\Psi(a + c-c', c, z)}{\Psi(a, c', z)}$$
for every $a > 0$ and $c' < c < 1$ is conjectured in Problem 7.2 in \cite{IK} and we will come back to this question in Remark \ref{Single} (d) below.} 
\end{Rem}

We next obtain the following reinforcement of Corollary \ref{B0} in the case $a\le 1.$ The argument relies on an appropriate expression for the cumulative Thorin measure of the beta prime distribution. 

\begin{Coro}
\label{B2}
For every $a \in (0,1]$ and $c' < c < 1,$ the function
$$z\; \mapsto\; \frac{\Psi(a, c, z)}{\Psi(a, c', z)}$$  
is {\em LCM}.
\end{Coro}

\proof
Taking the pointwise limit, it is enough to consider the case $a\in (0,1).$ The LCM property amounts to the ID property of the above random variable $\Y_{a,b,b'},$ where we have again set $b = 1-c$ and $b' = 1-c'.$ Recall from the self-decomposability of $\B'_{a,b'}$ that the random variable $\X^{a,b'}_x$ in the above proof is ID for every $x > 0.$ On the other hand, the ID property may not be preserved by mixtures. From the above decomposition, we see that the log-Laplace transform of $\Y_{a,b,b'}$ is given by
\begin{equation}
\label{LapY}
-\log \esp\lcr e^{-z \Y_{a,b,b'}}\rcr \; =\; \int_0^\infty (1 - e^{-zy}) (\mu_{a,b} - \mu_{a, b'})(dy)
\end{equation}
where $\mu_{a,x}$ is the L\'evy measure of $\B'_{a,x}$ for $x = b, b'.$ The density of $\mu_{a,x}$ on $(0,\infty)$ can be recovered from (1.4) of \cite{IK}, which implies from (\ref{LT}) and the L\'evy-Khintchine formula that (here and throughout, $\Psi'_z$ denotes the derivative of $\Psi$ w.r.t. to $z$)
$$\int_0^\infty y\,e^{-zy}\, \mu_{a,x} (dy) \; = \; -\frac{\Psi'_z(a,1-x,z)}{\Psi(a, 1-x,z)}\; = \; \frac{a\Psi(a+1, 2-x,z)}{\Psi(a,1-x,z)}\; =\; a\int_0^\infty\!\! e^{-zy} \,\esp[e^{-y \GG_{a,x}}]\,dy$$
where $\GG_{a,x}$ is a positive random variable with density
\begin{equation}
\label{AWK}
\frac{t^{x-1} e^{-t} \lva \Psi(a,1-x, t e^{\ii \pi})\rva^{-2}}{\Ga(a+1)\Ga(a+x)}\cdot
\end{equation}
In \cite{JRY}, this random variable is called the Thorin random variable associated to $\B'_{a,x},$ as mentioned in the introduction. By uniqueness of the Laplace transform, the density of $\mu_{a,x}(dy)$ on $(0,\infty)$ reads 
$$a\, y^{-1} \esp\lcr e^{-y \GG_{a,x}}\rcr \; = \; a \int_0^\infty \!\! e^{-yt}\, \pb [ \GG_{a,x} \le t] \, dt$$
for every $a,x > 0$ and we can deduce from (\ref{LapY}) that the LCM property is a consequence of the stochastic ordering
\begin{equation}
\label{STO}
\GG_{a,b}\; \prost \; \GG_{a,b'}
\end{equation}
for every $a \in (0,1)$ and $b' > b > 0.$\footnote{Throughout, for two positive random variables $X, Y,$ we will use the standard notation $X\prost Y$ to express the fact that $\pb[X\le x] \le \pb[Y\le x]$ for all $x \ge 0$.} Indeed, the ordering (\ref{STO}) implies that the measure $\mu_{a,b} - \mu_{a, b'}$ in (\ref{LapY}) is positive and integrates $1\wedge y$ and is hence a L\'evy measure. In Remark \ref{Single} below we will see that (\ref{STO}) can be derived directly, though not immediately, from (\ref{AWK}) for every $a\in(0,1),$ and we will also conjecture a stronger property. 

Anticipating the argument for Corollary \ref{Para2} below, we prefer using a general closed formula for finite Thorin measures which has been obtained in \cite{JRY}. Specifically, it follows from (51), (53) and Theorem 2.3(4) therein that for every $a \in(0,1)$ and $t,x >0,$ one has
$$\pb [ \GG_{a,x} \le t] \; = \; 1\, -\, \frac{1}{\pi a} \arctan \lpa \frac{\sin(\pi a)}{\cos(\pi a) + f_{a,x}(t) }\rpa\; =\; \frac{\sin(\pi a)}{\pi a} \int_0^{f_{a,x}(t)} \!\!\!\frac{du}{u^2 + 2 \cos (\pi a) u +1}$$
with
$$f_{a,x}(t) \; =\; \frac{\esp [(\G_x^{-1} - t^{-1})_+^{-a}]}{\esp[(t^{-1} - \G_x^{-1})_+^{-a}]}\cdot$$  
We next compute 
$$\esp [(\G_x^{-1} - t^{-1})_+^{-a}] \; = \; \frac{1}{\Ga(x)} \int_{t^{-1}}^\infty (y - t^{-1})^{-a} y^{-x-1} e^{-y^{-1}} dy\; = \; \frac{t^{a+x} e^{-t}}{\Ga(x)} \int_0^1 u^{-a} (1-u)^{a+x-1}\, e^{tu}\, du$$
by the change of variable $y = 1/t(1-u),$ and 
$$\esp [(t^{-1} - \G_x^{-1})_+^{-a}] \; = \; \frac{1}{\Ga(x)} \int_0^{t^{-1}} (t^{-1} - y)^{-a} y^{-x-1} e^{-y^{-1}} dy\; = \; \frac{t^{a+x} e^{-t}}{\Ga(x)} \int_0^\infty \!u^{-a} (1+u)^{a+x-1} e^{-tu}\, du$$
by the change of variable $y = 1/t(1+u).$ Therefore, 
\begin{equation}
\label{RatioInt}
f_{a,x}(t) \; =\; \frac{{\displaystyle \int_0^1 u^{-a} (1-u)^{a+x-1}\, e^{tu}\, du}}{{\displaystyle \int_0^\infty\! u^{-a} (1+u)^{a+x-1} e^{-tu} \,du}}
\end{equation}
and it is plain from this formula that $x\mapsto f_{a,x}(t)$ decreases on $(0,\infty),$ for every $a\in (0,1)$ and $t > 0.$ Putting everything together, we obtain the required stochastic ordering (\ref{STO}).

\endproof

\begin{Rem}
\label{Single}
{\em (a) The function $f_{a,x}(t)$ in the above proof has the alternative expression
$$f_{a,x}(t) \; = \;\frac{\Ga(a+x)\, \Phi(1-a, 1+x,t)}{\Ga(1+x)\, \Psi(1-a,1+x,t)}$$
where $\Phi (a,c,z)$ is Kummer's confluent hypergeometric function - see Formula 6.5.(1) in \cite{EMOT}. The function $f_{a,x}$ is well-defined on $(0,\infty)$ for $a= 1$ or more generally when $a = n$ a positive integer thanks to Formula 6.9.2.(36) in \cite{EMOT}, which yields $f_{n,x}\equiv (-1)^{n-1}.$ On the other hand, $f_{a,x}$ has a finite number of poles on $(0,\infty)$ for $a> 1$ not an integer, because $\Psi(1-a, 1+x,t)$ has then a finite number of zeroes which are different from those of $\Phi(1-a, 1+x, t)$ on $(0,\infty)$ - see Paragraph 6.16 in \cite{EMOT}.

It is clear from (\ref{RatioInt}) that $t\mapsto f_{a,x}(t)$ is an increasing bijection onto $(0,\infty).$ This property also follows from the computation of the Wronskian 
$$\Phi'_t(1-a, 1+x,t)\Psi(1-a,1+x,t)\, -\, \Phi(1-a, 1+x,t)\Psi'_t(1-a,1+x,t)\; =\;\frac{\Ga(1+x)}{\Ga(1-a)}\, t^{-1-x} e^t $$
which is given in Formula 6.7.(6) in \cite{EMOT} with the constant $K_{15}$ therein. This Wronskian also allows one to deduce from (\ref{AWK}) the formula
\begin{equation}
\label{Cumul}
\pb [ \GG_{a,x} \le t] \; = \; \frac{\sin(\pi a)}{\pi a} \int_0^{f_{a,x}(t)} \!\!\!\frac{du}{u^2 + 2 \cos (\pi a) u +1}
\end{equation}
for the cumulative Thorin measure of $\B'_{a,x}$ for every $a\in (0,1)$ and $x > 0,$ which was used in the above proof. Indeed, we have
\begin{eqnarray*}
\frac{t^{x-1} e^{-t} \lva \Psi(a,1-x, t e^{\ii \pi})\rva^{-2}}{\Ga(a+1)\Ga(a+x)} & = & \frac{\sin(\pi a)}{\pi a} \lva\frac{\Ga(1-a)\, t^x\, e^{-t}\,\Psi(1-a, 1+x,t)}{\Ga(a+x) \Psi(a,1-x, t e^{\ii \pi})}\rva^2 f_{a,x}'(t)\\
& = & \frac{\sin(\pi a)}{\pi a} \lva\frac{\Ga(1-a)\, \Psi(1-a, 1+x,t)}{\Ga(a+x)\, e^t\,\Psi(a+x,1+x, t e^{\ii \pi})}\rva^2 f_{a,x}'(t)\\
& = & \frac{\sin(\pi a)\,f_{a,x}'(t)}{\pi a \lva f_{a,x}(t) \, +\, e^{\ii \pi a}\rva^{2}}\\ 
& = & \frac{\sin(\pi a)\,f_{a,x}'(t)}{\pi a\, (f_{a,x}^2 (t) + 2 \cos(\pi a) f_{a,x}(t) + 1)}
\end{eqnarray*}
where we have used Formul\ae\, 6.5.(6) resp. 6.7.(7) in \cite{EMOT} for the second resp. third equality. This implies (\ref{Cumul}) by integration since $t\mapsto f_{a,x}(t)$ is an increasing bijection onto $(0,\infty).$ 

We mention that this formula (\ref{Cumul}) corrects (145) in \cite{JRY}, which was obtained after an erroneous application of the formula (127) therein. More precisely, it is easy to see that Formula (145) in \cite{JRY} amounts to the alternative formula 
$$f_{a,x}(t) \; =\; \frac{{\displaystyle \int_0^\infty u^{-a} (1+u)^{a+x-1}\, e^{-tu^{-1}}\, du}}{{\displaystyle \int_0^1\! u^{-a} (1-u)^{a+x-1} e^{tu^{-1}} \,du}}$$
which is not true since this would imply that $x\mapsto f_{a,x}(t)$ increases for all $a\in (0,1)$ and $t >0,$ contradicting \eqref{STO}.
\\

(b) The formula (\ref{Cumul}) also makes it possible to compute the cumulative Thorin measure of the Pareto random variable\footnote{also called {\em Lomax} in the literature} $\B'_{1,x}$ for every $x >0.$ To do so, observe first that by Formula 6.7.(7) in \cite{EMOT}, one has 
$$\frac{f_{1-\varepsilon,x}(t) -1}{\varepsilon\pi}\;\to\; -\lpa\ii\, +\, \frac{\Ga(1+x) e^{t + \ii\pi x}}{\pi}\, \Psi(1+x, 1+x, e^{\ii\pi} t)\rpa$$
as $\varepsilon\downarrow 0.$ On the other hand, the same formula 6.7.(7) in \cite{EMOT} shows that
\begin{eqnarray*}
\frac{\Ga(1+x) e^{t + \ii\pi x}}{\pi}\, \Psi(1+x, 1+x, e^{\ii\pi} t) & = & \lim_{a\downarrow 0} \frac{1}{\pi a} \lpa \frac{\Ga(1+x) e^{-\ii\pi a}}{\Ga(1+x-a)} \,\Psi(a,1+x,t)\, -\, \Phi(a,1+x,t)\rpa\\
& = & \frac{1}{\pi} \int_0^\infty y^{-1}\!\lpa(1+y)^x e^{-ty} - (1-y)_+^x e^{ty}\rpa dy\, -\, \ii.
\end{eqnarray*}
Considering now the integral of the Frullani type
$$g_x(t)\; =\;  \frac{1}{\pi} \int_0^\infty y^{-1}\!\lpa(1-y)_+^x e^{ty} - (1+y)^x e^{-ty}\rpa dy,$$
we see that $t\mapsto g_x(t)$ is an increasing bijection from $(0,\infty)$ onto $\rl$ and, from the above discussion, that $f_{1-\varepsilon, x} (t) = 1 +\varepsilon \pi g_x(t) + o(\varepsilon).$ The change of variable $u = 1 + \pi\varepsilon v$ in (\ref{Cumul}) implies then
$$\pb[\GG_{1-\varepsilon,x} \le t] \; = \; \frac{\sin (\pi\varepsilon)}{\pi^2\varepsilon} \int_{(\pi\varepsilon)^{-1}}^{g_x(t) +o(1)}\!\! \frac{dv}{v^2 + 1 + o(1)}$$
and, letting $\varepsilon\to 0,$ we finally obtain the relatively simple formula
$$\pb[\GG_{1,x} \le t] \; = \; \frac{1}{\pi} \int_{-\infty}^{g_x(t)}\!\! \frac{dv}{v^2+1}\cdot$$
This corrects the complicated formula (153) in \cite{JRY}, which is a consequence of  (145) therein and is hence also erroneous by the above discussion. Observe the identity $\GG_{1,x}\elaw g_x^{-1}(\C)$ where $\C$ is the standard Cauchy random variable. Combining this identity with (\ref{AWK}) in the case $a = 1$ and Formula 6.9.2.(21) in \cite{EMOT}, we also obtain the following formula for the modulus of the incomplete Gamma function on the cut
$$\frac{t^{-x-1} e^{t} \lva \Gamma(-x, t e^{\ii \pi})\rva^{-2}}{\Ga(x+1)}\; =\; \frac{g_x'(t)}{\pi (1 + g_x^2(t))}$$
for every $t,x > 0.$ Let us finally mention the convergence $\GG_{a,x}\claw \G_x$ as $a\downarrow 0,$ a direct consequence of (\ref{AWK}) and Scheff\'e's lemma.\\

(c) The so-called Askey-Wimp-Kerov distribution is a symmetric distribution on $\rl$ with density
$$w_c(t)\; =\; \frac{\lva D_{-c}(\ii t)\rva^{-2}}{\sqrt{2\pi}\, \Ga(1+c)}$$
where $D_{-c}$ is the classical parabolic cylinder function with parameter $c > -1.$ It was proved in Theorem 4 of \cite{AW} that $w_c(t)$  is the weight function on $\rl$ associated to an orthogonal family of generalized Hermite polynomials. An exponential representation of its Stieltjes transform is also given in Theorem (8.2.5) of \cite{Ke}, which was recently applied in \cite{HST} to show that the normal distribution is freely self-decomposable. It is worth mentioning that the Askey-Wimp-Kerov distribution is, for $c > 0,$ a symmetrization of the case $x=1/2$ in (\ref{AWK}). More precisely, it follows from Formula 6.9.2.(31) in \cite{EMOT} that $w_c(t)$ is the density of the random variable 
$${\bf W}_c\; =\; \varepsilon \sqrt{ 2\GG_{c/2,1/2}}$$ where $\varepsilon$ is a Rademacher random variable that is $\pb[\varepsilon = 1] = \pb[\varepsilon =\!-1] = 1/2,$ independent of $\GG_{c/2,1/2}.$ Observe that its extends at the limit $c =0,$ with ${\bf W}_0$ having a standard normal distribution and $\GG_{0,1/2} \elaw \G_{1/2}.$ \\

(d) Since $\B'_{a,b} \elaw\B'_{a',b} \times\B_{a,a'-a}$ for every $b > 0$ and $a' > a > 0,$ we have $\B'_{a'\!,b}\,\prost\, \B'_{a,b}$ and it is natural, in view of Corollary \ref{B1}, to conjecture that for every $c < 1$ and $a' > a > 0,$ the function
$$z\; \mapsto\; \frac{\Psi(a', c, z)}{\Psi(a, c, z)}$$  
is CM. Notice that if the latter is true, then by Corollary \ref{B1} the function
$$z\; \mapsto\; \frac{\Psi(a+\lbd, c, z)}{\Psi(a, c-\mu, z)}$$
is also CM for every $a,\lbd,\mu >0$ and $c < 1.$ As mentioned above in Remark 1, this is formulated as Problem 7.2. in \cite{IK} in the case $\lbd =\mu,$ and this problem is still open to the best of our knowledge. See also the subsequent Problem 7.3. in \cite{IK}. Observe from our above argument that the question would be solved if we could prove that $a'\pb [ \GG_{a',b} \le t] \ge a\pb [ \GG_{a,b} \le t]$ for every $b,t > 0$ and $a' > a > 0.$ \\

(e) We believe that the statement of Corollary \ref{B2} holds for every $a > 0$ and $c < 1$. Unfortunately, there does not seem to exist any tractable formula for the distribution function of $\GG_{a,x}$ which would help prove that $x\mapsto\pb[\GG_{a,x}\le t]$ is non-increasing also when $a > 1.$ Observe that by (\ref{AWK}) and the single intersection property - see Theorem 1.A.12 in \cite{SS}, this stochastic ordering would also be a consequence of 
$$\sharp\lacc t \in (0,\infty)\;\;\;\big/\;\;\;t^{c-c'}\lva\frac{\Psi(a,c, t e^{\ii \pi})}{\Psi(a,c', t e^{\ii \pi})} \rva^2\, =\, \frac{\Ga(1+a-c')}{\Ga(1+a-c)}\racc\; =\; 1$$
for every $a > 0$ and $c\neq c' <1,$ a fact which is supported by simulations but still eludes us.}
\end{Rem}

\section{Proof of Theorem A}

We begin with a general formula for the density of $\lbd \B'_{a,b} + \mu\B'_{c,d}$ in terms of the first Appell series $F_1$, which is more convenient for our purposes than the ones previously obtained in \cite{D,E} for the distribution function. This formula will also be used later on in Section 5.  We refer to \cite{AK} and to sections 5.7-14 in \cite{EMOT} for classic accounts on Appell series.

\begin{Propo}
\label{Main0}
For every $a,b,c,d,\lbd,\mu > 0,$ the density of $\lbd \B'_{a,b} + \mu\B'_{c,d}$ reads 
$$\frac{\Ga(c+d)\lbd^{-a}\mu^d\, x^{a+c-1}}{\Ga(c)\Ga(d)\, (x+\mu)^{c+d}}\, F_1\lpa a, a+b, c+d, a+c ; -x/\lbd, x/(x+\mu)\rpa$$
on $(0,\infty).$
\end{Propo}

\proof

The convolution formula and a change of variable show first that the density of $\lbd \B'_{a,b} + \mu\B'_{c,d}$ is  
\begin{equation}\label{referee1}
\frac{\Ga(a+b)\Ga(c+d)\lbd^b\mu^d\, x^{a+c-1}}{\Ga(a)\Ga(b)\Ga(c)\Ga(d)}\,\int_0^1\! \frac{y^{a-1} (1-y)^{c-1}}{(\lbd + xy)^{a+b}\, (x+\mu -xy)^{c+d}}\, dy.
\end{equation}
By Picard's integral formula - see Formul\ae\, 5.7.(6) and 5.8.2.(5) in \cite{EMOT}, we obtain the required expression as an Appell series of the first kind. 

\endproof

We next give a simpler expression in terms of the hypergeometric function for the density of $\B'_{a,b} + \B'_{a,b}$ which we will, throughout, denote by $\pab.$ This expression will be our starting point for the hypergeometric identities which will be used henceforth.
  
\begin{Propo}
\label{Main1}
For every $a,b > 0$ one has
$$\pab(x)\; =\; \frac{\Ga(a+b)^2\, x^{2a-1}}{\Ga(2a) \Ga(b)^2\, (x+1)^{a+b}}\; \pFq{2}{1}{a+b,,a}{a+1/2}{\frac{-x^2}{4(x+1)}}$$
on $(0,\infty).$
\end{Propo}

\proof

Letting $\lbd =\mu = 1$ and $(a,b) = (c,d)$ in \eqref{referee1}, the integral transforms into
$$\frac{2\,\Ga(a+b)^2\, x^{2a-1}}{\Ga(a)^2\Ga(b)^2}\,\int_0^{1/2}\! \frac{(y(1-y))^{a-1}}{(x+1 + x^2 y(1-y))^{a+b}}\, dy$$
and then, making the substitution $z = 4y(1-y),$ into
$$\frac{2^{1-2a}\,\Ga(a+b)^2\, x^{2a-1}}{\Ga(a)^2\Ga(b)^2\, (x+1)^{a+b}}\,\int_0^1\! \frac{z^{a-1}(1-z)^{-1/2}}{(1 + x^2z/4(x+1))^{a+b}}\, dz.$$
We can then conclude by Euler's integral formula and Legendre's multiplication formula - see respectively Formul\ae\, 2.1.3.(10) and 1.2.(15) in \cite{EMOT}.

\endproof

\begin{Rem}
{\em The above change of variable also gives the following expression for the density of the independent sum $\B_{a,b} +\B_{a,b},$ which reads
$$\frac{\Ga(a+b)^2\, x^{2a-1}}{\Ga(2a) \Ga(b)^2\, (1-x)^{1-b}}\; \pFq{2}{1}{1-b,,a}{a+1/2}{\frac{-x^2}{4(1-x)}}$$
on $(0,1)$ and is obtained on $(1,2)$ by the same formula with the switchings $a\leftrightarrow b$ and $x\leftrightarrow 2-x.$}

\end{Rem}
We now proceed to the proof of Theorem A. In the case $b = 1/2,$ the hypergeometric function in (\ref{Main1}) simplifies by the binomial theorem and gives the following formula, which could have been obtained directly 
\begin{equation}
\label{Dens12}
\varphi_{a,1/2}(x) \; =\; \frac{2 \Ga(a+1/2)\, x^{2a-1}}{\Ga(a) \sqrt{\pi}\, (x+1)^{1/2} (x+2)^{2a}}\cdot
\end{equation}
By (\ref{Dens12}), for every $s \in (-2a, 1/2)$ the Mellin transform reads
\begin{eqnarray*}
\esp\lcr \lpa \B'_{a,1/2} + \B'_{a,1/2}\rpa^s \rcr & = & \frac{2 \Ga(a+1/2)}{\Ga(a) \sqrt{\pi}}\, \int_0^\infty\frac{x^{2a +s -1}}{(x+1)^{1/2} (x+2)^{2a}}\; dx\\
& = & \frac{2 \Ga(a+1/2)}{\Ga(a) \sqrt{\pi}}\, \int_0^1\frac{u^{-1/2-s} (1-u)^{2a +s -1}}{(1+u)^{2a}}\; du\\
& = & \frac{2 \Ga(a+1/2)\Ga(2a)}{\Ga(a) \Ga(2a+1/2)}\, \pFq{2}{1}{2a,,1/2-s}{2a+1/2}{-1}\,\times\, \esp\lcr \lpa \B'_{2a,1/2}\rpa^s \rcr,
\end{eqnarray*}
where in the third equality we have used (\ref{Mell}) and Euler's integral formula. We next transform the hypergeometric expression into
\begin{eqnarray*}
\frac{2 \Ga(a+1/2)\Ga(2a)}{\Ga(a) \Ga(2a+1/2)}\, \pFq{2}{1}{2a,1/2-s}{2a+1/2}{-1} & = & \frac{2 \Ga(a+1/2)\Ga(2a)}{\Ga(a) \Ga(2a+1/2)}\, \pFq{2}{1}{1/2-s,, 2a}{2a+1/2}{-1}\\
& = & \frac{2 \Ga(a+1/2)}{\Ga(a) \sqrt{\pi}}\, \int_0^1\frac{u^{2a-1}(1+u)^s}{(1-u^2)^{1/2}}\; du\\
& = & \frac{\Ga(a+1/2)}{\Ga(a) \sqrt{\pi}}\, \int_0^1\frac{u^{a-1}(1+\sqrt{u})^s}{(1-u)^{1/2}}\; du\\ 
& = & \esp \lcr \lpa 1 + \sqrt{\B_{a, 1/2}}\rpa^s\rcr.
\end{eqnarray*}
Identifying the Mellin factors, we obtain the required identity.

\qed

\medskip

This result implies the following CM property for ratios of certain confluent hypergeometric functions, which we will state here in terms of Hermite functions. Recall from Section 10.2 in \cite{L} that the Hermite function $H_\nu$ with parameter $\nu\in\rl$ is a generalization of the classical Hermite polynomials $H_n,$ given by
\begin{equation}
\label{DH}
H_\nu (z)\; =\; 2^{\nu/2} e^{z^2/2} D_\nu (\sqrt{2} z)
\end{equation}
where $D_\nu$ stands for the parabolic cylinder function. In the case of a negative parameter one has also the integral representation
\begin{equation}
\label{Herman}
H_{-\nu}(z)\; =\; \frac{1}{\Ga(\nu)} \int_0^\infty e^{-t^2 -2tz}\, t^{\nu-1} dt\; =\; \frac{\Ga(\nu/2)}{2\Ga(\nu)} \; \esp\lcr e^{-2z\sqrt{\G_{\nu/2}}}\rcr
\end{equation}
for every $\nu > 0$ and $z\in\rl$ - see Formula (10.5.2) in \cite{L}.

\begin{Coro}
\label{Para1}
For every $\nu > 0,$ the function
$$z\; \mapsto\; \frac{\lpa H_{-\nu}\lpa \sqrt{z}\rpa\rpa^2}{H^{}_{-2\nu}(\sqrt{z})}$$  
is {\em CM}.
\end{Coro}

\proof

Setting $a = \nu/2$ and applying (\ref{DH}), Formula 6.9.2.(31) in \cite{EMOT} and Theorem A, we have
$$\frac{\lpa H^{}_{-\nu}\lpa \sqrt{z}\rpa\rpa^2}{H^{}_{-2\nu}(\sqrt{z})}\; =\; \frac{\lpa\Psi(a, 1/2, z)\rpa^2}{\Psi(2a, 1/2,z)}\; =\; \frac{\Ga(2a+1/2)\sqrt{\pi}\;}{\lpa\Ga(a+1/2)\rpa^2}\times \frac{\esp\lcr e^{-z (\B'_{2a,1/2}(1 + \sqrt{\B_{a,1/2}}))}\rcr}{ \esp\lcr e^{-z\B'_{2a,1/2}}\rcr}\cdot$$
By the self-decomposability of $\B'_{2a,1/2},$ we can conclude exactly as in Corollary \ref{B1} via the independent decomposition 
$$\B'_{2a,1/2}\, \times\lpa 1 + \sqrt{\B_{a,1/2}}\rpa \; \elaw\; \B'_{2a,1/2} \, +\, \X_a^x$$
for some positive random variable $\X_a^x,$ conditionnally on $\B_{a,1/2} = x.$
 
\endproof

\begin{Rem}
\label{Para1b}
{\em Taking the limits at 0 and $\infty$, we obtain the following sharp bounds for the Hermite functions on $(0,\infty)$:
$$1\; <\;  \frac{\lpa H_{-\nu}(z)\rpa^2}{H^{}_{-2\nu}(z)}\; <\; \frac{(\Ga(\nu/2))^2\,\Ga(2\nu)}{2\,(\Ga(\nu))^3}\cdot$$
We refer to Remark \ref{DuranDuran} below for other sharp bounds on $(0,\infty)$ for ratios of Hermite functions with negative parameters.}
\end{Rem}

In the case $\nu \in(0,2],$ the CM property of Corollary \ref{Para1} can be refined.

\begin{Coro}
\label{Para2}
For every $\nu \in (0,2]$ the function
$$z\; \mapsto\; \frac{\lpa H_{-\nu}\lpa \sqrt{z}\rpa\rpa^2}{H^{}_{-2\nu}(\sqrt{z})}$$  
is {\em LCM}.
\end{Coro}

\proof By stability of the LCM property at the limit $\nu\uparrow 2,$ it suffices to consider the case $\nu\in(0,2).$ We set again $a = \nu/2\in (0,1).$ Similarly as in Corollary \ref{B2}, we have
$$-\log\lpa\frac{\lpa H_{-\nu}\lpa \sqrt{z}\rpa\rpa^2}{H^{}_{-2\nu}(\sqrt{z})}\rpa\; =\; \nu\int_0^\infty\!\! \int_0^\infty (1-e^{-zy})\; e^{-zt} \lpa \pb\lcr\HH_{a,1} \le t\rcr - \pb\lcr\HH_{a,2} \le t\rcr \rpa dt\, dy$$
with
$$\pb[\HH_{a,i} \le t]\; =\; \frac{\sin(\pi a)}{\pi a} \int_0^{g_{a,i}(t)} \!\!\!\frac{du}{u^2 + 2 \cos (\pi a) u +1}$$
for $i=1,2,$ where
$$g_{a,1}(t) \; =\; \frac{\esp \lcr\lpa \G_{1/2}^{-1}(1+ \sqrt{\B_{a/2,1/2}}) - t^{-1}\rpa_+^{-a}\rcr}{\esp\lcr\lpa t^{-1} - \G_{1/2}^{-1}(1+ \sqrt{\B_{a/2,1/2}})\rpa_+^{-a}\rcr}\qquad\mbox{and}\qquad g_{a,2}(t) \; =\; \frac{\esp \lcr\lpa \G_{1/2}^{-1} - t^{-1}\rpa_+^{-a}\rcr}{\esp\lcr\lpa t^{-1} - \G_{1/2}^{-1}\rpa_+^{-a}\rcr}\cdot$$
As above, it is enough to show that
$$g_{a,1}(t) \; \ge \;g_{a,2}(t) \; =\; \frac{{\displaystyle \int_0^1 y^{-a} (1-y)^{a-1/2}\, e^{ty}\, dy}}{{\displaystyle \int_0^\infty\! y^{-a} (1+y)^{a-1/2} e^{-ty} \,dy}}\cdot$$
To compute $g_{a,1}(t),$ we first need an evaluation of the density $\phi_a(x)$ of $\G_{1/2}^{-1}\,(1+ \sqrt{\B_{a/2,1/2}})$ which, since that of $1+ \sqrt{\B_{a/2,1/2}}$ reads
$$\frac{2 \Ga((a+1)/2) \, (x-1)^{a-1}}{\Ga(a/2)\sqrt{\pi\, x(2-x)}}$$ 
on $(1,2),$ is
\begin{eqnarray}
\label{phia}
\phi_a(x) & = & \frac{2 \Ga((a+1)/2)}{\pi \Ga(a/2)\, x^{3/2}} \, \int_0^1 y^{a-1} (1-y)^{-1/2} e^{-(y+1)x^{-1}}\! dy\nonumber\\
& = & \frac{2 \Ga((a+1)/2) \Ga(a) \, e^{-x^{-1}}}{\sqrt{\pi} \Ga(a/2)\Ga(a+1/2)\, x^{3/2}}\, \Phi(a,a+1/2, -x^{-1}).
\end{eqnarray}
Observe in passing that we recover when $a\downarrow 0$ the expression $\phi_0(x) = \frac{e^{-x^{-1}}}{\sqrt{\pi x^3}}$ for the density of $\G_{1/2}^{-1}$. We next compute
$$\esp \lcr\lpa \G_{1/2}^{-1}\lpa 1+ \sqrt{\B_{a/2,1/2}}\rpa - t^{-1}\rpa_+^{-a}\rcr \; = \; \int_{t^{-1}}^\infty \lpa y - t^{-1} \rpa^{-2a} \phi_a(y) \, dy$$
which equals, after some simplifications,
$$\frac{2 \Ga((a+1)/2) \Ga(a) t^{a-1/2} e^{-t}}{\sqrt{\pi} \Ga(a/2)\Ga(a+1/2)}\int_0^1 y^{-a}(1-y)^{a-1/2}\, e^{t y}\, \Phi (a,a+1/2, t(y-1))\, dy.$$
Similarly, we see that
$$\esp \lcr\lpa t^{-1} - \G_{1/2}^{-1}\lpa 1+ \sqrt{\B_{a/2,1/2}}\rpa\rpa_+^{-a}\rcr \; = \; \int^{t^{-1}}_0 \lpa t^{-1} - y \rpa^{-2a} \phi_a(y) \, dy$$
equals
$$\frac{2 \Ga((a+1)/2) \Ga(a) t^{a-1/2} e^{-t}}{\sqrt{\pi} \Ga(a/2)\Ga(a+1/2)}\int_0^\infty y^{-a}(1+y)^{a-1/2}\, e^{-t y}\, \Phi (a,a+1/2, -t(y+1))\, dy.$$
Putting the two blocks together, we obtain
$$g_{a,1}(t) \; = \; \frac{{\displaystyle \int_0^1\! y^{-a} (1-y)^{a-1/2}\, e^{ty}\, \Phi (a,a+1/2, t(y-1))\,dy}}{{\displaystyle \int_0^\infty\!\! y^{-a} (1+y)^{a-1/2} e^{-ty} \,\Phi (a,a+1/2, -t(y+1))\,dy}}\cdot$$
Using the bounds
$$\Phi (a,a+1/2, t(y-1))\;\ge\;\Phi (a,a+1/2, -t)\;\ge\;\Phi (a,a+1/2, -t(y+1))\; > \; 0$$
for every $t,y\ge 0,$ which are obvious consequences of Formula 6.5.(1) in \cite{EMOT}, we finally deduce the required inequality
$$g_{a,1}(t) \; \ge \;  \frac{{\displaystyle \int_0^1\! y^{-a} (1-y)^{a-1/2}\, e^{ty}\,dy}}{{\displaystyle \int_0^\infty\!\! y^{-a} (1+y)^{a-1/2} e^{-ty} \,dy}}\cdot$$
\endproof

\begin{Rem}
{\em (a) The formula (\ref{phia}) for $\phi_a$ is easily seen to be valid for all $a > 0.$ This implies by multiplicative convolution and Theorem A the following alternative formula for $\varphi_{a,1/2}$:
$$\varphi_{a,1/2} (x)\; =\; \frac{2 \Ga(a+1/2)\, x^{2a-1}}{\sqrt{\pi} \Ga(a)\Ga(2a+1/2)}\int_0^\infty y^{2a-1/2}\, e^{-(x+1)y}\, \Phi (2a,2a+1/2, -y)\, dy.$$
Setting $ z = x+1$ and comparing with (\ref{Dens12}), we recover the Laplace-Mellin transform
$$\int_0^\infty y^{2a-1/2}\, e^{-zy}\, \Phi (2a,2a+1/2, -y)\, dy\; =\; \frac{\Ga(2a+1/2)}{\sqrt{z}\, (z+1)^{2a}}$$
which is a particular case of the second Formula 6.10.(5) in \cite{EMOT}. Notice that this Laplace-Mellin transform and the above computation of $\varphi_a$ also give an indirect proof of Theorem A.\\

(b) With the notation of Section 2, one has $\HH_{a,1} \elaw \GG_{a,1/2}$ and $\HH_{a,2} \elaw \GG_{2a, 1/2}.$ Hence, the above proof shows that $\GG_{2a,1/2} \prost \GG_{a,1/2}$ for every $a\in (0,1).$ In general, we believe that the mapping $a\mapsto \pb[\GG_{a,1/2} \le t]$ is non-increasing on $(0,\infty)$ for every $t>0.$\\

(c) Since 
$$\frac{\lpa H_{-\nu}\lpa \sqrt{z}\rpa\rpa^2}{H^{}_{-2\nu}(\sqrt{z})}\;\to\; 1$$ 
as $z\to\infty,$ the ID random variable associated to this quotient for $\nu\in (0,2]$ by Corollary \ref{Para2} is compound Poisson, with an underlying finite L\'evy measure having density
$$\nu\! \int_0^\infty\!\! e^{-zt} \lpa \pb\lcr\HH_{\nu/2,1} \le t\rcr - \pb\lcr\HH_{\nu/2,2} \le t\rcr \rpa dt.$$}
%whose completely monotonic character means that there is also an exponential representation
%$$\frac{\lpa H_{-\nu}\lpa \sqrt{z}\rpa\rpa^2}{H^{}_{-2\nu}(\sqrt{z})}\; =\; \frac{(\Ga(\nu/2))^2\,\Ga(2\nu)}{2\,(\Ga(\nu))^3}\, \exp-\lpa\varphi_\nu (z)\rpa$$ 
%for some bounded complete Bernstein function $\varphi_\nu.$}

\end{Rem}

\section{Proof of Theorem B}

We start with two alternative expressions for $\pab$ which rely on Pfaff's transformations - see the two Formul\ae\, 2.1.4.(22) in \cite{EMOT} - applied to (\ref{Main1}):
\begin{equation}
\label{Pfaff1}
\pab(x)\; = \; \frac{4^{a+b} \Ga(a+b)^2\, x^{2a-1}}{\Ga(2a) \Ga(b)^2\, (x+2)^{2(a+b)}}\; \pFq{2}{1}{a+b,,1/2}{a+1/2}{\lpa \frac{x}{x+2}\rpa^2}
\end{equation}
and
\begin{equation}
\label{Pfaff2}
\pab(x)\; =\; \frac{4^a\, \Ga(a+b)^2\, x^{2a-1}}{\Ga(2a) \Ga(b)^2\, (x+1)^b(x+2)^{2a}}\; \pFq{2}{1}{1/2-b,,a}{a+1/2}{\lpa \frac{x}{x+2}\rpa^2}.
\end{equation}
The first expression (\ref{Pfaff1}) will be used for the case $a+b = 1,$ whereas the second expression (\ref{Pfaff2}) is our starting point to evaluate the Mellin transform 
$$\cM_{a,b} (s) = \esp \lcr \lpa \B'_{a,b} + \B'_{a,b} \rpa^s\rcr.$$ 
One has 
\begin{eqnarray}
\label{Melle}
\cM_{a,b} (s) & = & \frac{4^a\, \Ga(a+b)^2}{\Ga(2a) \Ga(b)^2} \int_0^\infty \frac{x^{2a+s-1}}{(x+1)^b(x+2)^{2a}}\; \pFq{2}{1}{1/2-b,,a}{a+1/2}{\lpa \frac{x}{x+2}\rpa^2}\, dx\nonumber\\
& = & \frac{2^{2a+s}\, \Ga(a+b)^2}{\Ga(2a) \Ga(b)^2} \int_0^1 \frac{y^{2a+s-1}(1-y)^{b-s-1}}{(y+1)^b}\; \pFq{2}{1}{1/2-b,,a}{a+1/2}{y^2}\, dy.
\end{eqnarray}
From this expression, we can first characterize the definition strip of $\cM_{a,b}(s)$.

\begin{Propo}
\label{MellF}
One has 
$$\cM_{a,b} (s)\, <\, \infty \quad \Longleftrightarrow\quad s\in (-2a,b).$$
\end{Propo}

\proof

By Gauss' summation formula - see Formula 2.1.3.(14) in \cite{EMOT}, we know that
$$\pFq{2}{1}{1/2-b,,a}{a+1/2}{1}\; = \; \frac{\Ga(a+1/2)\Ga(b)}{\sqrt{\pi}\,\Ga(a+b)}\; <\; \infty$$
and it follows from (\ref{Melle}) that
$$\cM_{a,b} (s)\, <\, \infty \quad \Longleftrightarrow\quad \int_0^1 y^{2a+s-1}(1-y)^{b-s-1}\, dy\, <\, \infty \quad \Longleftrightarrow\quad s\in (-2a,b).$$ 
\endproof

We next give a general formula for $\cM_{a,b}(s)$ in terms of some ${}_3F_2(1),$ which will be essentially used in the case $a = 1/2,$ but has an independent interest. 

\begin{Propo}
\label{3F2}
For every $a,b > 0$ and $s\in (-2a,b)$, one has
$$\cM_{a,b} (s)\; =\; \frac{2^{2a+s}\, \Ga(a+b)^2\, \Ga(2b-s)\,\Ga(2a+s)}{\Ga(2a) \Ga(b)^2\Ga(2a+2b)}\,\pFq{3}{2}{a+s/2,,a+(s+1)/2,,1/2}{a+1/2,,a+b+1/2}{1}.$$
\end{Propo}

\proof

Applying Euler's transformation 2.1.4.(23) in \cite{EMOT} to the hypergeometric function inside (\ref{Melle}) we obtain, by Fubini's theorem, 
\begin{eqnarray*}
\cM_{a,b} (s) & = & \frac{2^{2a+s}\, \Ga(a+b)^2}{\Ga(2a) \Ga(b)^2} \int_0^1 y^{2a+s-1} (1-y)^{2b-s-1}\; \pFq{2}{1}{a+b,,1/2}{a+1/2}{y^2}\, dy\\
& = & \frac{2^{2a+s}\, \Ga(a+b)^2}{\Ga(2a) \Ga(b)^2} \sum_{n\ge 0} \frac{(a+b)_n (1/2)_n}{(a+1/2)_n\, n!} \lpa \int_0^1 y^{2a+2n+s-1} (1-y)^{2b-s-1} dy\rpa\\
& = & \frac{2^{2a+s}\, \Ga(a+b)^2\, \Ga(2b-s)}{\Ga(2a) \Ga(b)^2} \sum_{n\ge 0} \frac{(a+b)_n (1/2)_n \Ga(2a+2n +s)}{(a+1/2)_n\, \Ga(2a + 2b + 2n)\,n!}\\
& = & \frac{2^{2a+s}\, \Ga(a+b)^2\, \Ga(2b-s)\,\Ga(2a+s)}{\Ga(2a) \Ga(b)^2\Ga(2a+2b)}\, \sum_{n\ge 0} \frac{(a+s/2)_n (a+(s+1)/2)_n(1/2)_n}{(a+1/2)_n(a + b + 1/2)_n\,n!}\\
& = & \frac{2^{2a+s}\, \Ga(a+b)^2\, \Ga(2b-s)\,\Ga(2a+s)}{\Ga(2a) \Ga(b)^2\Ga(2a+2b)}\,\pFq{3}{2}{a+s/2,,a+(s+1)/2,,1/2}{a+1/2,,a+b+1/2}{1}.
\end{eqnarray*}
\endproof

\subsection{The case $a + b = 1$}  Setting $\varphi_a = \varphi_{a,1-a}$ for concision, we first apply Kummer's quadratic transformation - see Formula 2.11.(28) in \cite{EMOT} - to (\ref{Pfaff1}) and get the following simple expression for the density: 
\begin{eqnarray*}
\varphi_a(x) & = & \frac{4\, x^{2a-1}}{\Ga(2a) \Ga(1-a)^2\, (x+2)^{2}}\; \pFq{2}{1}{1,,1/2}{a+1/2}{\lpa \frac{x}{x+2}\rpa^2}\\
& = & \frac{2\, x^{2a-1}}{\Ga(2a) \Ga(1-a)^2\, (x+2)}\; \pFq{2}{1}{1,,a}{2a}{-x}.
\end{eqnarray*}
If we next suppose $b < 1/2$ viz. $a > 1/2,$ Euler's integral formula implies
$$\varphi_a(x)\; =\; \frac{2\, x^{2a-1}}{\Ga(2a-1) \Ga(1-a)^2\, (x+2)}\,\int_0^1 \frac{(1-t)^{2a-2}}{(1+xt)^{a}}\, dt,$$
and we will use this expression to give an appropriate expression of the Mellin transform, via Fubini's theorem. Setting $\cM_a(s) = \cM_{a,1-a}(s)$ for concision we have, for every $s\in (-2a, 1-a),$
\begin{eqnarray*}
\cM_a(s) & = & \frac{2}{\Ga(2a-1) \Ga(1-a)^2}\,\int_0^1 (1-t)^{2a-2}\lpa \int_0^\infty \frac{x^{2a+s-1}}{(x+2)(1+xt)^{a}}\, dx \rpa dt\\
& = & \frac{2^{2a+s}}{\Ga(2a-1) \Ga(1-a)^2}\,\int_0^1 (1-t)^{2a-2}\lpa \int_0^1 \frac{z^{-a-s}(1-z)^{2a-1-s}}{(2t + z(1-2t))^{a}}\, dz \rpa dt\\
& = & \frac{2^{a+s}\Ga(2a+s)\Ga(1-a-s)}{\Ga(1+a)\Ga(2a-1) \Ga(1-a)^2}\,\int_0^1 t^{-a}(1-t)^{2a-2}\, \pFq{2}{1}{1-a-s,,a}{1+a}{1-\frac{1}{2t}} \,dt.
\end{eqnarray*}
Applying Euler's integral formula to the hypergeometric function inside the integral leads then to
$$\cM_a(s)\, = \,\frac{2\, (2a-1)}{\Ga(a)\Ga(1-a)}\lpa \int_0^1\!\!\int_0^1 t^{-a}(1-t)^{2a-2}u^{a-1} \lpa 2 + u(t^{-1} -2) \rpa^{s+a-1} dtdu \rpa \times\, \esp\lcr \lpa \B'_{2a,1-a}\rpa^s \rcr$$  
and, identifying the Mellin factors, to 
$$\B'_{a,1-a} \, +\, \B'_{a,1-a} \; \elaw\; \B'_{2a,1-a}\,\times \lpa 2 \, + \, \B_{a,1}\,\times\lpa \B_{1-a,2a-1}^{-1} -2\rpa\rpa^{(a-1)}$$
where we use, here and throughout, the standard notation $X^{(t)}_{}$ for the size-bias of order $t$ of a positive random variable $X$ with $\esp [X^t] < \infty,$ that is
\begin{equation}
\label{Sizemore}
\esp\lcr f(X^{(t)}_{})\rcr\; =\; \frac{\esp \lcr X^t_{} f(X)\rcr}{\esp[X^t_{}]}
\end{equation}
for every positive continuous function $f$ on $(0,\infty).$ Putting the blocks together, we are reduced to show the following strange identity in law.

\begin{Lemm}
\label{BetaStr}
For every $a\in (1/2,1),$ one has
$$\lpa 2 \, + \, \B_{a,1}\,\times\lpa \B_{1-a,2a-1}^{-1} -2\rpa\rpa^{(a-1)}\elaw\; 1 \, +\, \sqrt{\B_{a,1/2}\times\B_{1-a,a-1/2}^{-1}}\,\cdot$$
\end{Lemm}

\proof The argument is an evaluation of the densities of the r.v. $(2 +  \B_{a,1}\times (\B_{1-a,2a-1}^{-1} -2)$ and $(1 + (\B_{a,1/2}\times\B_{1-a,a-1/2}^{-1})^{1/2}),$ which we will respectively denote by $f_a$ and $g_a.$ To evaluate $g_a,$ we first compute the density of $\B_{a,1/2}\times\B_{1-a,a-1/2}^{-1},$ which reads
$$\frac{(a-1/2)\, x^{a-1}}{\Ga(a)\Ga(1-a)}\, \int_0^{1\wedge x^{-1}}\!\!\! (1-y)^{a-3/2} (1-xy)^{-1/2}\, dy.$$
and transforms into
$$\frac{x^{a-1}}{\Ga(a)\Ga(1-a)}\;  \pFq{2}{1}{1/2,,1}{a+1/2}{x}$$  
for $x < 1,$ and into
$$\frac{(2a-1)\, x^{a-2}}{\Ga(a)\Ga(1-a)}\,  \pFq{2}{1}{3/2-a,,1}{3/2}{x^{-1}}\; =\; \frac{(2a-1)\, (x-1)^{a-1}}{\Ga(a)\Ga(1-a)\, x}\,  \pFq{2}{1}{1/2,,a}{3/2}{x^{-1}}.$$
for $x > 1.$ Changing the variable, we then obtain
\begin{equation}
\label{12}
g_a(x)\; =\; \frac{2\,(x-1)^{2a-1}}{\Ga(a)\Ga(1-a)}\;\pFq{2}{1}{1/2,,1}{a+1/2}{(x-1)^2}
\end{equation}
on $(1,2)$ and
\begin{equation}
\label{2oo}
g_a(x)\; =\; \frac{2(2a-1) (x^2 -2x)^{a-1}}{\Ga(a)\Ga(1-a)\,(x-1)}\;\pFq{2}{1}{1/2,,a}{3/2}{(x-1)^{-2}}
\end{equation}
on $(2,\infty).$ We next proceed to the evaluation of $f_a$, which is more involved. The starting formula is the following expression for the density of $(\B_{1-a,2a-1}^{-1} -2):$
$$\frac{\Ga(a)}{\Ga(1-a)\Ga(2a-1)}\, (x+2)^{-a} (x+1)^{2a-2}, \qquad x > -1.$$
We deduce that the density of $\B_{a,1}\times(\B_{1-a,2a-1}^{-1} -2)$ equals 
$$\frac{\Ga(1+a)\, \vert x\vert^{a-1}}{\Ga(1-a)\Ga(2a-1)}\,\int_{\vert x\vert}^1 y^{-a}(2-y)^{-a} (1-y)^{2a-2} dy$$
on $(-1,0),$ which transforms into
$$\frac{\Ga(1+a)\, \vert x\vert^{a-1} (x+1)^{2a-1}}{2 \Ga(1-a)\Ga(2a-1)}\,\int_0^1 z^{a-3/2}(1-z(x+1)^2)^{-a}dz,$$
and then into
$$\frac{\Ga(1+a)\, \vert x\vert^{a-1} (x+1)^{2a-1}}{\Ga(1-a)\Ga(2a)}\,  \pFq{2}{1}{a,,a-1/2}{a+1/2}{(x+1)^2}.$$
Applying Euler's formula 2.1.4.(23) in \cite{EMOT} and translating the variable, we finally get
\begin{equation}
\label{12b}
f_a(x)\; =\; \frac{\Ga(1+a)\, x^{1-a}\, (x-1)^{2a-1}}{\Ga(1-a)\Ga(2a)}\,  \pFq{2}{1}{1/2,,1}{a+1/2}{(x-1)^2}
\end{equation}
on $(1,2).$ For $x > 0,$ the density of $\B_{a,1}\times(\B_{1-a,2a-1}^{-1} -2)$ equals
$$\frac{\Ga(1+a)\, x^{a-1}}{\Ga(1-a)\Ga(2a-1)}\,\int_x^\infty y^{-a}(2+y)^{-a} (1+y)^{2a-2} dy,$$
which transforms into
$$\frac{\Ga(1+a)\, x^{a-2}}{\Ga(1-a)\Ga(2a-1)}\,\int_0^1 \lpa 1+ 2yx^{-1}\rpa^{-a}\lpa 1+yx^{-1}\rpa^{2a-2}dy,$$
and then, by Picard's integral formula and Appell's transform - see Formula 5.8.2.(5) resp. 5.11.(4) in \cite{EMOT}, into
$$\frac{\Ga(1+a)\, (x+1)^{2a-2}(x+2)^{1-a}}{\Ga(1-a)\Ga(2a-1)\, x}\; F_1\lpa 1, a, 2-2a, 2 ; -2x^{-1}, -(x+1)^{-1}\rpa.$$
For $x\in(2,\infty),$ this implies
\begin{eqnarray}
\label{2oob}
f_a(x) & = & \frac{\Ga(1+a)\, x^{1-a} (x-1)^{2a-2}}{\Ga(1-a)\Ga(2a-1)\, (x-2)}\; F_1\lpa 1, a, 2-2a, 2 ; -2(x-2)^{-1}, -(x-1)^{-1}\rpa\nonumber \\
& = & \frac{\Ga(1+a)\, x^{-a} (x-1)^{2a-1}}{\Ga(1-a)\Ga(2a-1) (x-2)}\; F_1\lpa 1, a, a, 2 ; x^{-1}, (2-x)^{-1}\rpa \nonumber\\
& = & \frac{\Ga(1+a)\, (x-1)^{2a-1}(x-2)^{a-1}}{\Ga(1-a)\Ga(2a-1)}\,\int_0^1 ((x-1)^2 -y^2)^{-a} dy \nonumber\\
& = & \frac{(2a-1) \Ga(1+a)\,x^{1-a} (x^2-2x)^{a-1}}{\Ga(1-a)\Ga(2a)\, (x-1)}\; \pFq{2}{1}{1/2,,a}{3/2}{(x-1)^{-2}} 
\end{eqnarray}
where in the first equality we have used Appell's transform 5.11.(3) in \cite{EMOT}. Now, comparing (\ref{12}) to (\ref{12b}) resp. (\ref{2oo}) to (\ref{2oob}) shows that
$$g_a(x) \; = \; \frac{2\Ga(2a)}{a\,\Ga(a)^2} \; x^{a-1}f_a(x)$$
on $(1,2)\cup(2,\infty).$ This completes the proof by (\ref{Sizemore}).

\endproof

\begin{Rem}{\em In the case $b > 1/2$ viz. $a < 1/2,$ Euler's formula leads to
$$\varphi_a(x)\; =\; \frac{2\, x^{2a-1}}{\Ga(a)^2 \Ga(1-a)^2\, (x+2)}\,\int_0^1 \frac{(t(1-t))^{a-1}}{(1+xt)^{a}}\, dt.$$
Unfortunately, we could not derive from this expression of the density anything concrete on $\cM_a(s)$ which could lead to an explicit identity in law for $\B'_{a,1-a} + \B'_{a,1-a}$.}
\end{Rem}

\subsection{The case $a = 1/2$} Specifying Proposition \ref{3F2} to $a=1/2$, we have 
$$\cM_{a,b} (s)\; = \; \frac{2^{1+s}\, \Ga(b+1/2)^2\, \Ga(2b-s)\,\Ga(1+s)}{\Ga(b)^2\Ga(2b+1)}\;\pFq{3}{2}{1/2,,1+s/2,,(1+s)/2}{b+1,,1}{1}$$
for $s\in (-2a,b).$ By Thomae's relation - see 7.4.4.1 in \cite{BMP} or (2.2) in \cite{KR}, we deduce
$$\cM_{a,b} (s)\; = \;\frac{2^{1+s}\, \Ga(b+1/2)^2\, \Ga(2b-s)\,\Ga(1+s)\,\Ga(b-s)}{\sqrt{\pi}\, \Ga(b)^2\,\Ga(2b+1)\,\Ga(b+1/2-s)}\;\pFq{3}{2}{b-s/2,,b+(1-s)/2,,1/2}{b+1,,b+1/2-s}{1}.$$
The key-point is now that since $b+1/2-s = 2(b-s/2) +1/2 -(b+1) +1,$ we can apply the  quadratic transformation 7.4.1.11 in \cite{BMP} and get some ${}_3F_2(-1)$ instead of the ${}_3F_2(1)$ for $s\in[0,b):$ 
\begin{eqnarray*}
\cM_{a,b} (s) & = & \frac{2^{1+2b}\, \Ga(b+1/2)^2\, \Ga(2b-s)\,\Ga(1+s)\,\Ga(b-s)}{\sqrt{\pi}\, \Ga(b)^2\,\Ga(2b+1)\,\Ga(b+1/2-s)}\;\pFq{3}{2}{b-s,,b+1/2,,2b-s}{b+1,,b+1/2-s}{-1}.
\end{eqnarray*}
Observe that the hypergeometric series converges absolutely because $2b-s < 1.$ We now recall Pochhammer's integral formula
$$\pFq{3}{2}{\lbd,\mu,d}{\nu,e}{-z}\; =\; \frac{\Ga(\nu)\Ga(e)}{\Ga(\mu)\Ga(\nu -\mu)\Ga(d)\Ga(e-d)}\,\int_0^1\!\!\int_0^1 \frac{t^{\mu-1}(1-t)^{\nu-\mu-1} u^{d-1}  (1-u)^{e-d-1}}{(1+utz)^{\lbd}}\,dt du,$$
which is valid for every $\mu,d,z >0$ and $\nu > \mu, e > d,$ and is a consequence of Formula 2.1.3.(10) and Formul\ae\, 5.5.2.(5) and 5.6.(1) (with $p=q = 1$) in \cite{EMOT}. Applying this formula transforms the above expression for $\cM_{a,b}(s)$ into
$$\frac{2^{1+2b}\, \Ga(b+1/2)\Ga(1+s)\,\Ga(b-s)}{\pi\, \Ga(b)^2\,\Ga(2b+1)\,\Ga(1/2-b)}\int_0^1\!\!\int_0^1 t^{b-1/2}(1-t)^{-1/2} u^{b-1}  (1-u)^{-1/2-b}\,  (t+u^{-1})^{s-b}\,dt du,$$
which is finite for $s\in(-2a,b),$ and which we recognize as
$$\esp\lcr \lpa \B'_{1,b}\rpa^s\rcr\,\times\, \esp\lcr \lpa \lpa \B_{b+1/2,1/2} +\B_{b,1/2-b}^{-1}\rpa^{(-b)}\rpa^s\rcr.$$
Identifying the Mellin factors, we finally obtain
$$\B'_{1/2,b}\, +\, \B'_{1/2,b}\; \elaw\; \B'_{1,b}\,\times\, \lpa \B_{b+1/2,1/2} +\B_{b,1/2-b}^{-1}\rpa^{(-b)},$$
and we are reduced to prove the following identity in law.

\begin{Lemm}
\label{BetaStrb}
For every $b\in (0,1/2),$ one has
$$\lpa \B_{b+1/2,1/2} +\B_{b,1/2-b}^{-1}\rpa^{(-b)}\elaw\; 1 \, +\, \sqrt{\B_{1/2,1/2}\times\B_{1/2,1/2-b}^{-1}}\,\cdot$$
\end{Lemm}

\proof

As above, the proof relies on an evaluation of the densities of the r.v. $\B_{b+1/2,1/2} +\B_{b,1/2-b}^{-1}$ and $1 + (\B_{1/2,1/2}\times\B_{1/2,1/2-b}^{-1})^{1/2},$ which we will denote by $f_b$ and $g_b$ respectively. The computation of $g_b$ is done as in Lemma \ref{BetaStr} and, skipping details, we obtain
\begin{equation}
\label{b12}
g_b(x)\; =\; \frac{2\,\Ga(b+1/2)}{\sqrt{\pi}\,\Ga(b)}\;\pFq{2}{1}{1/2,,b+1/2}{1}{(x-1)^2}
\end{equation}
on $(1,2)$ and
\begin{equation}
\label{b2oo}
g_b(x)\; =\; \frac{2\,\Ga(b+1/2) \,(x -1)^{-2b-1}}{\Ga(b)\Ga(1+b)\Ga(1/2-b)}\;\pFq{2}{1}{b+1/2,,b+1/2}{b+1}{(x-1)^{-2}}
\end{equation}
on $(2,\infty).$ The evaluation of $f_b$ on $(1,\infty)$ starts with the formula
$$f_b(x)\; =\; \frac{b\, \cos(\pi b)}{\pi}\, \int_0^{1\wedge (x-1)}\!\!\! t^{b-1/2} (1-t)^{-1/2} (x-1-t)^{-1/2-b}(x-t)^{-1/2} dt$$ 
which is obtained by multiplicative convolution. On $(1,2),$ we have
\begin{eqnarray}
\label{b12b}
f_b(x) & = & \frac{b\,\cos(\pi b)}{\pi}\, \int_0^1 t^{b-1/2} (1-t)^{-1/2-b} (1-(x-1)t)^{-1/2}(x-(x-1)t)^{-1/2}\, dt\nonumber\\
& = & b\,x^{-1/2}\, F_1 \lpa b+1/2, 1/2, 1/2, 1; x-1, 1- x^{-1}\rpa\nonumber\\
& = & b\, x^b\; \pFq{2}{1}{1/2,,b+1/2}{1}{(x-1)^2},
\end{eqnarray}
where in the last equality we have used the reduction Formula 5.10.(1) in \cite{EMOT}. On $(2,\infty),$ the same argument shows 
\begin{eqnarray}
\label{b2oob}
f_b(x) & = & \frac{b\,(x-1)^{-1/2-b}}{\Ga(b+1/2)\,\Ga(1/2-b)\,\sqrt{x}}\, \int_0^1 t^{b-1/2} (1-t)^{-1/2} (1-(x-1)^{-1}t)^{-1/2-b}(1-x^{-1}t)^{-1/2}\, dt\nonumber\\
& = & \frac{\sqrt{\pi} \,(x-1)^{-1/2-b}}{\Ga(b)\,\Ga(1/2-b)\,\sqrt{x}}\, F_1 \lpa b+1/2, b+1/2, 1/2, b+1; (x-1)^{-1}, x^{-1}\rpa\nonumber\\
& = & \frac{\sqrt{\pi} \,(x-1)^{-2b-1}\, x^b}{\Ga(b)\,\Ga(1/2-b)}\, \pFq{2}{1}{b+1/2,,b+1/2}{b+1}{(x-1)^{-2}}.
\end{eqnarray} 
Comparing (\ref{b12}) to (\ref{b12b}) resp. (\ref{b2oo}) to (\ref{b2oob}) implies
$$g_b(x)\; =\; \frac{2\, \Ga(b+1/2)}{\sqrt{\pi}\, \Ga(b+1)} \; x^{-b}f_b(x)$$
on $(1,2)\cup(2,\infty),$ and concludes the proof as in Lemma \ref{BetaStr}.
\endproof

A direct consequence of Theorem B is the following result, extending Corollary \ref{Para1}.

\begin{Coro}
\label{CMpart}
For every $a > 0$ and $c\in[1/2,1],$ the function
$$z\;\mapsto \; \frac{\Psi^2(a,c,z)}{\Psi(2a,c,z)}$$
is {\em CM} if $a= 1/2, c=1/2,$ or $a=c.$
\end{Coro}

In Corollary \ref{b>1}, we will see that this CM property is however not true anymore for $c<0.$ On the other hand, we believe that it holds for any $a >0$ and $c\in [0,1]$ - see Conjecture \ref{CMcj} below. In the case $c\in (1/2,1),$ this would be a consequence of the following extension of Theorem B.
 
\begin{Conj} 
\label{CjMain}
For every $ a > 0$ and $b\in (0,1/2),$ one has
$$\B'_{a,b} \, + \, \B'_{a,b} \; \elaw \; \B'_{2a, b} \,\times \lpa 1 \, +\, \sqrt{\frac{\B_{a,1/2}}{\B_{b,1/2-b}}}\rpa.$$
\end{Conj}

As can be seen from Proposition \ref{3F2}, this conjecture is equivalent to the following non-conventional integral representation of a certain ${}_3F_2(1),$ to be compared with the aforementioned Pochhammer formula: 
$$\pFq{3}{2}{a+s/2,,a+(s+1)/2,,1/2}{a+1/2,,a+b+1/2}{1}\; =\; \frac{\Ga(b-s)\,\Ga(a+1/2)\,\Ga(a+b+1/2)}{\Ga(a)\, \Ga(1/2-b)\,\Ga(a+b)\, \Ga(b-s/2)\,\Ga(b+(1-s)/2)}\qquad\qquad$$
\begin{flushright}
$\times\;{\displaystyle \int_0^1\!\!\int_0^1 u^{a-1}\,(1-u)^{-1/2}\, v^{b-1}\,(1-v)^{-1/2-b} \lpa 1+\sqrt{\frac{u}{v}}\rpa^s du\, dv,}$
\end{flushright}

\noindent
or alternatively
$$\pFq{3}{2}{a+s/2,,a+(s+1)/2,,1/2}{a+1/2,,a+b+1/2}{1}\; =\; \frac{4\,\Ga(b-s)\,\Ga(a+1/2)\,\Ga(a+b+1/2)}{\Ga(a)\, \Ga(1/2-b)\,\Ga(a+b)\, \Ga(b-s/2)\,\Ga(b+(1-s)/2)}\qquad\qquad$$
\begin{flushright}
$\times\;{\displaystyle \int_0^1\!\!\int_0^1 u^{2a-1}\,(1-u^2)^{-1/2}\, v^{2b-s-1}\,(1-v^2)^{-1/2-b}\, (u + v)^s du\, dv}$
\end{flushright}

\medskip

\noindent
for every $a > 0, b \in (0,1/2)$ and $s\in (-2a,b).$ The non-conventional character of this representation comes from the square root inside the integral. In the cases $a = 1/2$ and $a = 1-b,$ we were able to show the representation essentially with the help of Thomae's relations on ${}_3F_2(1)$ and the reduction of Appell series of the first kind, our main tools in the above proof. Unfortunately, we do not know as yet how to extend such arguments to the general case. There is a huge list of classic transformations on ${}_3F_2(z)$ displayed in Section 7.4. of \cite{BMP}, and in the case $a=1/2$ we have made a crucial use of Formula 7.4.1.11. therein. Our feeling is that these classic transformations might not be enough to tackle Conjecture \ref{CjMain} in general. In this respect, let us refer to \cite{KR} for a recent and non-classical transformation on some ${}_3F_2(1),$ which however does not seem to be very helpful for our problem. See also Remark \ref{ConjHyp} below for another integral representation of some ${}_2F_1(z)$ which would give an answer to Conjecture \ref{CjMain} in the case $a+b =1/2.$

\section{Miscellaneous}

\subsection{The case $a+ b=1/2$} In this paragraph we prove another identity in law for the convolution of $\B'_{a,b}$ with itself in the case $a+b = 1/2,$ which is interesting in its own right. Notice that this identity would seem to imply the validity of Conjecture \ref{CjMain} in this case as well - see Remark \ref{ConjHyp} below. Unfortunately, this fact still eludes us.

\begin{Propo}
\label{ab12}
If $a+b =1/2,$ one has
$$\B'_{a,b} \, +\, \B'_{a,b}\;\elaw\; \B'_{2a,2b}\,\times\lpa \B_{a,1/2} \, +\, \B_{b,1/2}^{-1}\rpa.$$
\end{Propo}

\proof
Specifying Proposition \ref{Mell} to $a+b =1/2,$ we get
$$\cM_{a,b} (s)\; =\; \frac{2^{2a+s} \pi\, \Ga(2b-s)\,\Ga(2a+s)}{\Ga(2a) \Ga(b)}\,\pFq{3}{2}{a+s/2,,a+(s+1)/2,,1/2}{a+1/2,,1}{1}.$$
As in the case $a=1/2,$ we next obtain a formula with a ${}_3 F_2(-1)$ instead of the  ${}_3 F_2(1):$
\begin{eqnarray*}
\cM_{a,b} (s) & = & \frac{2^{2a+s} \sqrt{\pi}\, \Ga(2b-s)\,\Ga(b-s)\,\Ga(2a+s)}{\Ga(2a) \Ga(b)\,\Ga(b+1/2 -s)}\,\pFq{3}{2}{-s/2,,(1-s)/2,,1/2}{a+1/2,,b+1/2-s}{1}\\
& = & \frac{2^{1-2b} \sqrt{\pi}\, \Ga(2b-s)\,\Ga(b-s)\,\Ga(2a+s)}{\Ga(2a) \Ga(b)\,\Ga(b+1/2 -s)}\,\pFq{3}{2}{-s,,a,,b-s}{a+1/2,,b+1/2-s}{-1}.
\end{eqnarray*}
Applying finally the Pochhammer integral formula and the Legendre multiplication formula, we transform this expression into
$$\frac{\Ga(a+1/2)\,\Ga(b+1/2)\,\Ga(2a+s)\,\Ga(2b-s)}{\pi\,\Ga(a)\, \Ga(b)\,\Ga(2a)\,\Ga(2b)}\int_0^1\!\!\int_0^1 t^{a-1}(1-t)^{-1/2} u^{b-1}  (1-u)^{-1/2}\,  (t+u^{-1})^s\,dt du,$$
which is
$$\esp\lcr \lpa \B'_{2a,2b}\rpa^s\rcr\,\times\, \esp\lcr \lpa \B_{a,1/2} +\B_{b,1/2}^{-1}\rpa^s\rcr.$$ 
\endproof

As a by-product, we obtain in passing the following result on the ratio of a Macdonald function of order zero and the exponential integral, which also has a doubling character. 
 
\begin{Coro}
\label{Expi}
The function
$$z\;\mapsto\; \frac{\lpa K^{}_0(z) \rpa^2}{E^{}_1(2z)}$$
is {\em CM}.
\end{Coro}

\proof

Reasoning as in Corollary \ref{B1}, it follows from (\ref{LT}), Proposition \ref{ab12} and the self-decomposability of $\B'_{2a,2b}$ that the function
$$z\; \mapsto\; \frac{\lpa \Psi (a,1/2-a,2z)\rpa^2}{\Psi(2a,2a,2z)}\; =\; \frac{e^{2z} (2z)^{1-2a} (K_{a-1/2} (z))^2}{\pi\, \Psi(2a,2a,2z)}$$
is CM for every $a\in (0,1/2),$ where the equality is a consequence of Formul\ae\, 6.9.(2) and 7.2.2.(13) in \cite{EMOT}. On the other hand, a change of variable shows that
$$\Psi(2a,2a,2z)\; =\; \frac{e^{2z}}{\Ga(2a)}\, \int_{2z}^\infty e^{-t}\lpa\frac{t}{2z} -1\rpa^{2a-1}\frac{dt}{t}$$
and the right-hand side converges to $e^{2z}E_1(2z)$ pointwise as $a\to 1/2.$ This completes the proof by the stability of the CM property under pointwise convergence.

\endproof

\begin{Rem}
\label{ConjHyp}
{\em It is natural to ask if Proposition \ref{ab12} implies Conjecture \ref{CjMain}. It is clear from (\ref{base}) and (\ref{base1}) that this implication amounts to the identity 
$$\B_{a,1/2} \, +\, \B_{b,1/2}^{-1} \; \elaw \;\B_{b,b}^{-1}\,\times\, \lpa 1 \, +\, \sqrt{\frac{\B_{a,1/2}}{\B_{b,a}}}\rpa$$
for $a+b = 1/2.$ It can be shown that the density of the r.v. on the left-hand side equals
$$\frac{\Ga(b+1/2)\,(x-1)^{a-1/2}}{\sqrt{\pi}\, \Ga(b)\, x}\; \pFq{2}{1}{1/2,,a}{a+1/2}{\frac{4(x-1)}{x^2}}$$
on $(1,2)\cup(2,\infty)$ and is remarkably given by the same function on both intervals - a function which diverges at 2. The density of the r.v. on the right-hand side is a priori more complicated and equals
$$\frac{2^{1-2a}\, \Ga(b+1/2)\,(x-1)^{a-1/2}}{\Ga(a)\, \Ga(b)^2\, x^{2b}} \int_0^1 y^{-2b}(1-y)^{b-1} (1-(x-1)y)^{-b} \,\pFq{2}{1}{a,,a}{a+1/2}{((x-1)y)^2}\, dy$$
on $(1,2).$ In particular, the two densities have the same asymptotics as $x\downarrow 1$ and $x\uparrow 2,$ respectively given by
$$\frac{\Ga(b+1/2)\,(x-1)^{a-1/2}}{\sqrt{\pi}\, \Ga(b)}\qquad\quad\mbox{and}\quad\qquad \frac{\Ga(a+1/2)\,\Ga(b+1/2)}{\pi\,\Ga(a)\,\Ga(b)}\, \log\lpa \frac{1}{2-x}\rpa.$$
However, we are currently unable to prove that the two densities are equal on $(1,2).$ Notice that this is equivalent to the integral formula
$$\frac{\sqrt{\pi}}{\Ga(a)\, \Ga(b)} \int_0^1 y^{2a-1}(1-y)^{b-1} (1-zy)^{-b} \pFq{2}{1}{a,a}{a+1/2}{(zy)^2} dy\; =\; \lpa\!\frac{z+1}{2}\!\rpa^{2b}\!\pFq{2}{1}{1/2,a}{a+1/2}{\frac{4z}{(z+1)^2}}$$
for $z\in(0,1),$ which does not seem to follow from e.g. Erd\'elyi's fractional integral - see Formula 2.4.(3) in \cite{EMOT} - combined with some quadratic transform. We leave this question open to future investigations.}
\end{Rem}

\subsection{The case $b >1$} In this paragraph we show that an identity similar to Conjecture \ref{CjMain} cannot hold when $b > 1$ - see Corollary \ref{b>1} below. This fact is a consequence of the following stochastic ordering property, which has an independent interest.  

\begin{Propo}
\label{Stoo}
For every $a > 0,$ one has
$$\B'_{a,b}\, +\,\B'_{a,b}\;\,\prost \;\,\B'_{2a,b}\quad \Longleftrightarrow\quad b \le 1.$$
\end{Propo}

\proof

Setting $\tpab$ for the density of $\B'_{2a,b}$ and recalling the formula (\ref{Main1}) for the density $\pab$ of $\B'_{a,b} +\B'_{a,b},$ we have
\begin{eqnarray*}
\frac{\tpab(x)}{\pab(x)} & = & \frac{\Ga(a+b)^2\, (x+1)^{a}}{\Ga(2a+b) \Ga(b)}\; \pFq{2}{1}{a+b,,a}{a+1/2}{\frac{-x^2}{4(x+1)}}\\
& = &  \frac{\Ga(a+b)^2\, (z+1)^{2a}}{\Ga(2a+b) \Ga(b)}\; \pFq{2}{1}{1/2-b,,a}{a+1/2}{z^2}
\end{eqnarray*}
with $z = x/(x+2),$ where in the second equality we have used Pfaff's transformation. The limits of this function as $x\to 0$ and $x\to \infty$ are respectively
$$\lambda_{a,b}\; =\; \frac{\Ga(a+b)^2}{\Ga(2a+b)\,\Ga(b)}\qquad\mbox{and}\qquad \Lambda_{a,b}\; =\; \frac{4^a\,\Ga(a+1/2)\,\Ga(a+b)}{\sqrt{\pi}\, \Ga(2a+b)}\cdot$$
By the strict log-convexity of the Gamma function, one has $\lbd_{a,b} < 1$ and the function $b\mapsto \Lambda_{a,b}$ decreases on $(0,\infty).$ Hence, since $\Lambda_{a,1} = 1$ by the Legendre multiplication formula, one also has $\Lambda_{a,b} < 1$ if $b > 1$ and this finishes the proof of the only if part: if $b > 1,$ there exists by continuity of $\pab$ and $\tpab$ some $x_\ast > 0$ such that $\pab(x)< \tpab(x)$ for every $x \ge x_\ast,$ which implies
$$\pb\lcr \B'_{a,b} +\B'_{a,b} \ge x_\ast\rcr\; < \; \pb\lcr \B'_{2a,b} \ge x_\ast\rcr$$
and contradicts $\B'_{a,b} +\B'_{a,b}\,\prost\, \B'_{2a,b}.$ For the if part, since $\lbd_{a,b} < 1$ it is enough to show by the single intersection property that if $b\le 1,$ one has
\begin{equation}
\label{Card}
\sharp\{ x >0, \; \tpab(x) =\pab(x)\}\; =\; 1.
\end{equation}
If $b\le 1/2,$ it is plain from the above that the function 
$$z\;\mapsto\; \psi_{a,b}(z)\; =\; \frac{\tpab(x)}{\pab(x)}\; =\;\frac{\Ga(a+b)^2\, (z+1)^{2a}}{\Ga(2a+b) \Ga(b)}\; \pFq{2}{1}{1/2-b,,a}{a+1/2}{z^2} $$ increases on $(0,1)$ from $\lbd_{a,b} < 1$ to $\Lambda_{a,b} > 1,$ so that (\ref{Card}) holds. If $b\in(1/2,1],$ since $\Lambda_{a,b} \ge 1$ it is clear that (\ref{Card}) will hold as soon as $z\mapsto\psi_{a,b}(z)$ increases and then decreases. We next compute the derivative $\psi'_{a,b}(z),$ which equals
$$\frac{2a\,\Ga(a+b)^2\, (z+1)^{2a-1}}{\Ga(2a+b) \Ga(b)}\,\lpa \pFq{2}{1}{1/2-b,,a}{a+1/2}{z^2}\, -\, \frac{(b-1/2)\, z(z+1)}{a+1/2}\, \pFq{2}{1}{3/2-b,,a+1}{a+3/2}{z^2}\rpa$$  
on $(0,1)$ and we see that
$$\psi'_{a,b}(0)\, =\, \frac{2a\,\Ga(a+b)^2}{\Ga(2a+b) \Ga(b)}\, > \, 0\qquad\mbox{and}\qquad \psi'_{a,b}(1)\, =\, -\infty.$$      
Finally, since
$$z\;\mapsto\; \pFq{2}{1}{1/2-b,,a}{a+1/2}{z^2}\qquad\mbox{and}\qquad z\;\mapsto\; -\lpa\pFq{2}{1}{3/2-b,,a+1}{a+3/2}{z^2}\rpa$$
clearly decrease on $(0,1)$ when $b\in(1/2,1],$ we deduce that $z\mapsto\psi'_{a,b}(z)$ changes its sign once on $(0,1)$, which implies (\ref{Card}) and completes the proof.

\endproof

\begin{Coro}
\label{b>1}
For every $a> 0$ and $c<0,$ the function
$$z\;\mapsto \; \frac{\Psi^2(a,c,z)}{\Psi(2a,c,z)}$$
is not {\em CM}. In particular, for every $a > 0$ and $b>1$ there does not exist any positive random variable $\X$ such that
$$\B'_{a,b} \, +\, \B'_{a,b}\;\elaw\; \B'_{2a,b}\,\times\lpa 1 \, +\, \X\rpa.$$
\end{Coro}

\proof

Set $b = 1-c > 1.$ If the above function were CM, then by (\ref{LT}) and Bernstein's theorem, we would have an additive factorization
$$\B'_{a,b}\, +\,\B'_{a,b}\;\elaw \;\B'_{2a,b}\, +\, \Y$$
for some positive r.v. $\Y,$ which is excluded by the only if part of Proposition \ref{Stoo}. The second negative statement follows as above from the self-decomposability of $\B'_{2a,b},$ by contraposition.
 
\endproof

\subsection{Monotonicity properties} In this paragraph we establish a monotonicity property for the ratio of $\Psi$ functions appearing in Conjecture \ref{CjMain}, and we state a related open problem. We also give a simple proof of the Tur\'an inequality for the $\Psi$ function and the parabolic cylinder function recently obtained in \cite{BI, K}, together with a fractional extension. All these results rely on the following consequence of a classical lemma by Biernaki and Krzy{\. z} \cite{BK}, which might be well-known but we could not find any reference.

\begin{Lemm}
\label{LTMon}
Let $f,g$ be probability densities defined on an interval $(a,b)\subseteq (0,\infty)$ which are positive, continuous, and such that the ratio $x\mapsto f(x)/g(x)$ is non-decreasing on $(a,b).$ Then, the ratios
$$z\;\mapsto \;\frac{{\displaystyle \int_0^\infty\!\! e^{-zx} f(x)\, dx}}{{\displaystyle \int_0^\infty\!\! e^{-zx} g(x)\, dx}}\qquad\quad\mbox{and}\qquad\quad z\;\mapsto \;\frac{{\displaystyle \int_0^\infty\!\! \frac{f(x)}{(1+xz)^\mu}\, dx}}{{\displaystyle \int_0^\infty\!\! \frac{g(x)}{(1+xz)^\mu}\, dx}}$$
are non-increasing on $(0,\infty)$ for every $\mu > 0.$ 
\end{Lemm}

\proof

We start with the ratio of Laplace transforms. Suppose first $b < \infty$ and fix some $n\in\NN^*.$ Consider the variate $u = e^{-z(b-a)/n}$ and, for every $i = 1,\ldots, n,$ set $x_i = a + i(b-a)/n, a_i = f(x_i)$ and $b_i = g(x_i).$ Since the sequence $\{a_i/b_i\}$ is non-decreasing, it follows from Biernaki and Krzy{\. z}'s lemma that the function
$$u\;\mapsto\; \frac{{\displaystyle \sum_{i=1}^n a_i u^i}}{{\displaystyle \sum_{i=1}^n b_i u^i}}$$
is non-decreasing on $(0,1),$ so that the function
$$z\;\mapsto\;\frac{{\displaystyle \sum_{i=1}^n e^{-z x_i} f(x_i)}}{{\displaystyle \sum_{i=1}^n e^{-z x_i} g(x_i)}}$$
is non-increasing on $(0,\infty).$ We can then conclude by taking the pointwise limit of the Riemann sums as $n\to\infty.$ The case $b=\infty$ is analogous and we omit details. To handle the ratio of generalized Stieltjes transforms, it suffices to observe by a change of variable that 
$$\frac{{\displaystyle \int_0^\infty\!\! \frac{f(x)}{(1+xz)^\mu}\, dx}}{{\displaystyle \int_0^\infty\!\! \frac{g(x)}{(1+xz)^\mu}\, dx}}\; =\; \frac{{\displaystyle \int_0^\infty\!\! e^{-z^{-1}u} \lpa u^{\mu -1} \int_0^\infty f(x)\, e^{-xu} dx\rpa du }}{{\displaystyle \int_0^\infty\!\! e^{-z^{-1} u}  \lpa u^{\mu -1} \int_0^\infty g(x)\, e^{-xu} dx\rpa du}}\cdot$$
Hence, the monotonicity of the left-hand side can be deduced from that of the ratios of Laplace transforms. More precisely, the ratios of functions between brackets inside the integrals on the right-hand side is non-increasing in $u$, so that the ratio on the left-hand side is non-decreasing in $z^{-1}$ and non-increasing in $z$.

\endproof

\begin{Rem}
{\em For the ratio of Laplace transforms, it is easy to check that the monotonicity property extends on $(c, \infty)$ with $c\in [-\infty,0)$, provided both Laplace transforms are finite on $(c,\infty)$.}
\end{Rem}

\begin{Coro}
\label{b<1/2}
For every $a> 0$ and $c\in [1/2,1],$ the function
$$z\;\mapsto \; \frac{\Psi^2(a,c,z)}{\Psi(2a,c,z)}$$
decreases on $(0,\infty).$
\end{Coro}

\proof
Set $b = 1 - c\in [0,1/2]$ and suppose first $b > 0.$ Using (\ref{LT}) and the notation of the proof of Proposition \ref{Stoo}, we have
$$\frac{\Psi^2(a,c,z)}{\Psi(2a,c,z)}\; =\; \frac{\Ga(b)\Ga(2a+b)}{(\Ga(a+b))^2}\,\times\lpa \frac{{\displaystyle \int_0^\infty\!\! e^{-zx} \pab(x)\, dx}}{{\displaystyle \int_0^\infty\!\! e^{-zx} \tpab(x)\, dx}}\rpa.$$
We have seen during the proof of Proposition \ref{Stoo} that the function
$$x\;\mapsto\; \frac{\pab(x)}{\tpab(x)}$$
is non-increasing on $(0,\infty)$ for every $a > 0$ and $b\in (0,1/2],$ so that we can conclude by Lemma \ref{LTMon} that 
$$z\;\mapsto \; \frac{\Psi^2(a,c,z)}{\Psi(2a,c,z)}$$
is non-increasing. This property remains true for $b=0$ by pointwise limit. Finally, it is clear by analyticity that the above function cannot be locally constant, so that it actually decreases.

\endproof

In view of Corollaries \ref{b>1} and \ref{b<1/2}, it is natural to state the following open problem, which might be challenging. Recall that in the case $c\in [1/2,1]$ the open if part of this statement is, by self-decomposability, a consequence of Conjecture \ref{CjMain}.

\begin{Conj}
\label{CMcj}
For every $a> 0,$ the function
$$z\;\mapsto \; \frac{\Psi^2(a,c,z)}{\Psi(2a,c,z)}$$
is {\em CM} if and only if $c\in [0,1].$
\end{Conj}

As another, independent, consequence of Lemma \ref{LTMon}, we obtain the following extension of a recent result  - see Theorem 3.1 in \cite{K} - on ratios of parabolic cylinder functions, with a very simple argument. This generalization has the same flavour as the main results of \cite{MS} on Kummer and Gauss hypergeometric functions.
 
\begin{Coro}
\label{Hermit}
For every $\nu, c > 0,$ the function
$$z\;\mapsto \; \frac{(H_{-\nu-c}(z))^2}{H_{-\nu}(z) H_{-\nu -2c}(z)}\; =\; \frac{(D_{-\nu-c}(\sqrt{2} z))^2}{D_{-\nu}(\sqrt{2} z) D_{-\nu -2c}(\sqrt{2} z)}$$
decreases on $\rl.$
\end{Coro}

\proof
By the same analyticity argument as in Corollary \ref{b<1/2}, it is enough to show that the ratio is non-increasing. Set $a = (\nu +c)/2$ and $t = c/2.$ By (\ref{Herman}), we have the representation
$$\frac{(H_{-\nu-c}(z))^2}{H_{-\nu}(z) H_{-\nu -2c}(z)}\; =\; \frac{\Ga(a)^2\, \Ga(2a-2t)\,\Ga(2a+2t)}{\Ga(2a)^2\, \Ga(a-t)\, \Ga(a+t)}\,\times \lpa \frac{\esp\lcr e^{-2z(\sqrt{\G_a} + \sqrt{\G_a})}\rcr}{\esp\lcr e^{-2z(\sqrt{\G_{a-t}} + \sqrt{\G_{a+t}})}\rcr}\rpa$$
for all $z\in \rl.$ On the other hand, the density of $\!\sqrt{\G_a} + \!\sqrt{\G_a}$ is easily computed as
\begin{eqnarray*}
h_a(x) & = & \frac{4 x^{4a-1} e^{-x^2}}{\Ga(a)^2} \int_0^1 (y(1-y))^{2a-1} e^{2 x^2 y(1-y)}\, dy\\
& = & \frac{2 x^{4a-1} e^{-x^2}}{\Ga(a)^2} \int_0^1 \lpa \frac{z}{4}\rpa^{2a-1}\! (1-z)^{-1/2}\, e^{x^2z/2}\, dz
\end{eqnarray*}
where in the equality we have used the substitution $z=4y(1-y)$ as in the beginning of Section 3. Similarly, the density of $\!\sqrt{\G_{a-t}} +\!\sqrt{\G_{a+t}}$ is computed as 
$$h_{a,t}(x) \; =\; \frac{2 x^{4a-1} e^{-x^2}}{\Ga(a)^2} \int_0^1 \lpa \frac{z}{4}\rpa^{2a-1}\! (1-z)^{-1/2} \lpa u_t(z) + u_t(z)^{-1}\rpa\,  e^{x^2z/2}\, dz$$
with the notation
$$u_t(z)\; =\; \lpa \frac{1+\sqrt{1-z}}{1-\sqrt{1-z}}\rpa^{2t}.$$
Since $z\mapsto u_t(z)$ is decreasing from $(0,1)$ to $(1,\infty)$ and since $u\mapsto u + u^{-1}$ is increasing on $(1,\infty),$ we deduce from the above formul\ae\, and Lemma \ref{LTMon} that the ratio $x\mapsto h_a(x)/h_{a,t}(x)$ is non-decreasing on $(0,\infty),$ and we can conclude by the remark made after Lemma \ref{LTMon}. 

\endproof

\begin{Rem}
\label{DuranDuran}
{\em As in Corollary 3.2 of \cite{K}, taking the limits at $\pm\infty$ implies the following sharp bounds
$$1\; <\; \frac{(D_{-\nu-c}(z))^2}{D_{-\nu}(z)\, D_{-\nu -2c}(z)} \; <\; \frac{\Ga(\nu)\,\Ga(\nu +2c)}{(\Ga(\nu +c))^2}$$ 
which are valid for every $\nu, c > 0$ and $z\in\rl.$ These bounds can be viewed as ``fractional'' Tur\'an inequalities for the parabolic cylinder function. See the recent paper \cite{Seg1} and the references therein for related results on ratios of parabolic cylinder functions.}
\end{Rem}

As a last consequence of Lemma \ref{LTMon}, we obtain the following fractional Tur\'an inequality for the $\Psi$ function, which generalizes Theorem 2 in \cite{BI}. 

\begin{Coro}
\label{Trico}
For every $a, \lbd > 0$ and $c <1,$ the function
$$z\;\mapsto \; \frac{\Psi(a, c-2\lbd, z) \Psi(a+2\lbd, c, z)}{(\Psi (a+\lbd,c-\lbd, z))^2}$$
decreases on $(0,\infty).$ In particular, one has the following sharp inequalities
$$ 1\; <\;  \frac{\Psi(a, c-2\lbd, z) \Psi(a+2\lbd, c, z)}{(\Psi (a+\lbd,c-\lbd, z))^2}\; < \; \frac{\Ga(1-c)\Ga(1-c +2\lbd)}{\Ga(1-c+\lbd)^2}\cdot$$
\end{Coro}

\proof

Again, it suffices to prove that the function is non-increasing. Observe indeed that the sharp inequalities are obtained in taking the limits at 0 and $\infty,$ as in Remark \ref{DuranDuran}. Setting $b = 1-c >0,$ we know by (\ref{LT}) that
$$\frac{\Psi(a, c-2\lbd, z) \Psi(a+2\lbd, c, z)}{(\Psi (a+\lbd,c-\lbd, z))^2}\; =\; \frac{\Ga(1-c)\Ga(1-c +2\lbd)}{\Ga(1-c+\lbd)^2}\times \frac{\esp\lcr e^{-z(\B'_{a,b+2\lbd} + \B'_{a+2\lbd,b})}\rcr}{\esp\lcr e^{-z(\B'_{a+\lbd,b+\lbd} + \B'_{a+\lbd,b+\lbd})}\rcr}\cdot$$
By Proposition \ref{Main0}, the densities of $\B'_{a,b+2\lbd} + \B'_{a+2\lbd,b}$ and $\B'_{a+\lbd,b+\lbd} + \B'_{a+\lbd,b+\lbd}$ read respectively
$$\frac{\Ga(a+b+2\lbd)\, x^{2a+2\lbd-1}}{\Ga(b)\Ga(a+2\lbd)\, (x+1)^{a+b+2\lbd}}\, F_1\lpa a, a+b+2\lbd, a+b+2\lbd, 2a+2\lbd ; -x, x/(x+1)\rpa$$
and
$$\frac{\Ga(a+b+2\lbd)\, x^{2a+2\lbd-1}}{\Ga(b+\lbd)\Ga(a+\lbd)\, (x+1)^{a+b+2\lbd}}\, F_1\lpa a+\lbd, a+b+2\lbd, a+b+2\lbd, 2a+2\lbd ; -x, x/(x+1)\rpa$$
on $(0,\infty).$ Therefore, by Lemma \ref{LTMon}, it suffices to prove that the ratio
$$x\;\mapsto\; \frac{F_1\lpa a, a+b+2\lbd, a+b+2\lbd, 2a+2\lbd ; -x, x/(x+1)\rpa}{F_1\lpa a +\lbd, a+b+2\lbd, a+b+2\lbd, 2a+2\lbd ; -x, x/(x+1)\rpa}$$
is non-decreasing on $(0,\infty).$ Setting $X(x) = x^2/4(x+1)$ for concision, we next observe by Picard's integral formula for Appell series, and some simplifications analogous to those of the proofs of Proposition \ref{Main1} and Corollary \ref{Hermit}, that the above ratio is proportional to 
$$\frac{{\displaystyle \int_0^1\! z^{a+\lbd -1}(1-z)^{-1/2}\,(1 + z X(x))^{-(a+b+2\lbd)}\, \lpa u_{\lbd/2} (z) + u_{\lbd/2}(z)^{-1}\rpa\, dz }}{{\displaystyle \int_0^1\! z^{a+\lbd -1}(1-z)^{-1/2}\,(1 + z X(x))^{-(a+b+2\lbd)}\, dz }}$$
which is a non-decreasing function of $x$ by the increasing character of $x\mapsto X(x),$ the decreasing character of $z\mapsto u_{\lbd/2} (z) + u_{\lbd/2}(z)^{-1},$ and Lemma \ref{LTMon}. 

\endproof

\begin{Rem}{\em We refer to \cite{BPS, TS} for some other, non-fractional, refinements of the Tur\'an inequalities obtained in \cite{BI}. See in particular Corollaire 2 in \cite{TS} and Theorem 1 in \cite{BPS} for different improvements of Theorem 2 in \cite{BI}.}
\end{Rem}

\subsection{Multiplicative convolution} The random variable 
$$\lpa 1\,+\,\B'_{a,b}\rpa\,\times\,\lpa 1\,+\,\B'_{c,d}\rpa\; -\; 1 \; =\; \B'_{a,b}\;+\;\B'_{c,d}\; +\; \B'_{a,b}\B'_{c,d}$$
is tighty related to the additive convolution $\B'_{a,b}+\B'_{c,d}$ that we have considered in the present paper, through the addition of the non-independent term $\B'_{a,b}\B'_{c,d}.$ It is worth mentioning the following identity in law for this random variable, since it is valid without restrictions on the parameters. In the case $(a,b) = (c,d),$ observe in particular the multiplicative factor $\B'_{2a, b}$ which appears in the same position as in Theorems A and B.

\begin{Propo}
\label{Free}
For every $a,b,c,d > 0$ such that $b < c+ d,$ one has
$$\lpa 1\,+\,\B'_{a,b}\rpa\,\times\,\lpa 1\,+\,\B'_{c,d}\rpa\;-\; 1\;\elaw\; \B'_{a+c,d}\,\times \lpa 1\, +\, \B_{a,c}\,\B'_{c+d-b, b}\rpa.$$ 
\end{Propo}

\proof

By (\ref{base}), the identity amounts to 
$$\frac{1}{\B_{b,a}\B_{d,c}}\, -\, 1\;\elaw\;\B'_{a+c,d}\,\times \lpa 1\, +\, \B_{a,c}\,\B'_{c+d-b, b}\rpa.$$ 
By Pochhammer's integral formula and some simplifications, we have
$$\esp\lcr\lpa \frac{1}{\B_{b,a}\B_{d,c}}\, -\, 1\rpa^s\rcr\; =\; \frac{\Ga(b-s)\Ga(d-s)\Ga(a+b)\Ga(c+d)}{\Ga(b)\Ga(d)\Ga(a+b-s)\Ga(c+d-s)}\;\pFq{3}{2}{-s,,b-s,,d-s}{a+b-s,,c+d-s}{1}$$
for $s \in (-a-c, b\wedge d).$ Applying Thomae's relation 7.4.4.2. in \cite{BMP}, we deduce
$$\esp\lcr\lpa \frac{1}{\B_{b,a}\B_{d,c}}\, -\, 1\rpa^s\rcr\; =\; \frac{\Ga(a+c+s)\Ga(d-s)\Ga(a+b)\Ga(c+d)}{\Ga(b+d)\Ga(b)\Ga(d)\Ga(a+c+d)}\;\pFq{3}{2}{a,,c+d-b,,a+c+s}{a+c,,a+c+d}{1}.$$
Using now Pochhammer's formula in the reverse direction, the right-hand side transforms into
$$\esp\lcr\lpa \B'_{a+c,d}\rpa^{-s}\rcr\,\times\, \frac{\Ga(a+c)\Ga(c+d)}{\Ga(a)\Ga(b)\Ga(c)\Ga(c+d-b)}\,\times\,\qquad\qquad\qquad\qquad\qquad\qquad\qquad\qquad\qquad\qquad\qquad\qquad$$
$$\qquad\qquad\qquad\qquad\qquad\qquad\qquad\int_0^1\!\int_0^1\! t^{c+d-b-1}(1-t)^{a+b-1}u^{a-1}(1-u)^{c-1}(1-ut)^{-b-d-s}\,dtdu$$
and we finally get, as in the proof of Theorem B, the intermediate identity
$$\lpa 1\,+\,\B'_{a,b}\rpa\,\times\,\lpa 1\,+\,\B'_{c,d}\rpa\;-\; 1\;\elaw\;\B'_{a+c,d}\;\times\;\lpa\frac{1}{1 - \B_{a,c}\B_{c+d-b,a+b}}\rpa^{(a+c)}.$$
We are hence reduced to show that
$$ 1\; +\;\B_{a,c}\,\B'_{c+d-b, b}\;\elaw\; \lpa\frac{1}{1 - \B_{a,c}\B_{c+d-b,a+b}}\rpa^{(a+c)}$$
which, as in Lemmas \ref{BetaStr} and \ref{BetaStrb}, will be obtained by an evaluation of the respective densities on $(1,\infty).$ The argument is here simpler since we do not need here to perform a separate evaluation on $(1,2)$ and on $(2,\infty).$ We show the equivalent identity 
\begin{equation}
\label{Hypergeo}
\frac{1}{1\, +\,\B_{a,c}\,\B'_{c+d-b, b}}\;\elaw\; \lpa 1 - \B_{a,c}\B_{c+d-b,a+b}\rpa^{(-a-c)},
\end{equation}
which, as will be observed in Remark \ref{Mult} (b) below, is a particular case of a general identity in law derived in \cite{LP}. Let us however give our own argument for the sake of completeness. The density of $1 - \B_{a,c}\B_{c+d-b,a+b}$ is standardly evaluated as
$$\frac{\Ga(a+c)\,\Ga(a+c+d)}{\Ga(a)\,\Ga(c+d-b)\,\Ga(a+b+c)}\; x^{a+b+c-1}(1-x)^{a-1}\pFq{2}{1}{a+b,,a+b-d}{a+b+c}{x},$$
and we see that the density of $(1 - \B_{a,c}\B_{c+d-b,a+b})^{(-a-c)}$ is proportional to 
$$x^{b-1}(1-x)^{a-1}\pFq{2}{1}{a+b,,a+b-d}{a+b+c}{x}.$$
On the other hand, the density of $(1+\B_{a,c}\,\B'_{c+d-b, b})^{-1}$ equals
$$\frac{\Ga(a+b)\,\Ga(a+c)\,\Ga(c+d)}{\Ga(a)\,\Ga(b)\,\Ga(c+d-b)\,\Ga(a+b+c)}\; x^{b-1}(1-x)^{-b-1}\pFq{2}{1}{a+b,,c+d}{a+b+c}{\frac{x}{x-1}},$$
which by Pfaff's transformation is also proportional to
$$x^{b-1}(1-x)^{a-1}\pFq{2}{1}{a+b,,a+b-d}{a+b+c}{x}.$$
\endproof

\begin{Rem}
\label{Mult}
{\em (a) The condition $b < c+d$ is not a restriction on the parameters: in the complementary situation $d < a+b,$ one has the factorization
$$\lpa 1\,+\,\B'_{a,b}\rpa\,\times\,\lpa 1\,+\,\B'_{c,d}\rpa\;-\;1\;\elaw\; \B'_{a+c,b}\,\times \lpa 1\, +\, \B_{a,c}\,\B'_{a+b-d, d}\rpa.$$

(b) The identity (\ref{Hypergeo}) is an instance of a general result of \cite{LP} for beta-hypergeometric distributions. With the notation of the introduction in \cite{LP}, one has indeed
\begin{eqnarray*}(1 - \B_{a,c}\B_{c+d-b,a+b})^{(-a-c)} & \sim  & {\rm BH} (b,a,a+b, a+b-d,a+b+c)\\
& \sim &  {\rm BH} (b,c+d-b,c+d,c,a+b+c)\;\sim\; {\rm BH} (Mv)
\end{eqnarray*}
with $v = (a,c,a+b-d,0,a+b),$ whereas $\B_{a,c}\sim {\rm BH} (v)$ and $\B'_{c+d-b,b}\sim \beta^{(2)}_{\Pi v}.$ Therefore, (\ref{Hypergeo}) is a consequence of (5) in \cite{LP}. It is worth mentioning that all size-biases of $1 - \B_{b,a}\B_{d,c}$ are beta-hypergeometric random variables, but that the converse is not true - see Theorems 2.1 and 2.2 in \cite{LP}. Observe finally that Thomae's relations play a crucial role in the proof of (5) in \cite{LP}.}
\end{Rem}

\subsection{The case of Mill's ratio}

The function $r : \rl \mapsto (0,\infty)$ defined by
$$r(x)\; =\; e^{x^2/2}\int_x^\infty e^{-t^2/2} dt$$
is known in the literature as Mill's ratio of the Gaussian law, whereas its reciprocal $1/r(x)$ is called the Gaussian hazard rate. This function, which appears countlessly in the literature, is related to the case $a=b=1/2$ of the present paper by the representations
\begin{equation}
\label{Mill}
r(\sqrt{2}\, x) \; = \; \sqrt{2}\, H_{-1}(x)\; =\; \sqrt{\frac{\pi}{2}} \, \esp\lcr e^{-2x\sqrt{\G_{1/2}}}\rcr, \qquad x\in\rl
\end{equation}
and
\begin{equation}
\label{Mill1}
r(\sqrt{2x})\; = \; \frac{1}{\sqrt{2}}\, \Psi(1/2,1/2,x)\; = \; \sqrt{\frac{\pi}{2}}\, \esp\lcr e^{-x\B'_{1/2,1/2}}\rcr, \qquad x\ge 0.
\end{equation}
The second equality in (\ref{Mill}) implies that $x\mapsto r(x)$ is LCM on $(0,\infty)$ by Corollary \ref{Para2}, and that it is up to normalization the MGF of the r.v. $-2\sqrt{\G_{1/2}},$ so that $x\mapsto r(x)$ is also log-convex on $\rl.$ The second equality in (\ref{Mill1}) entails that $x\mapsto r(\sqrt{x})$ is HCM on $(0,\infty)$ as the Laplace transform of a generalized Gamma convolution - see Theorem 5.4.1. in \cite{Bond}, and by Property v) p.68 in \cite{Bond} that the function $x\mapsto r(e^x)$ is log-concave on $\rl,$ or equivalently that the function $x\mapsto xr'(x)/r(x)$ is non-increasing on $(0,\infty).$ The above discussion gives hence a short proof of parts (a) and (b) of Theorem 2.5 in \cite{B}. Here is a simple improvement of parts (c) and (d) of this theorem.

\begin{Propo}
\label{Barr}
One has the following monotonicity properties.

{\em (a)} On $(0,\infty),$ the function $x\mapsto x^\a r(x)$ is decreasing for $\a = 0,$ increasing for $\a = 1$, and increasing then decreasing for $\a\in (0,1).$\

{\em (b)} On $(0,\infty),$ the function $x\mapsto x^\a r'(x)$ is increasing for $\a = 0,$ decreasing for $\a = 2$, and decreasing then increasing for $\a\in (0,2).$
\end{Propo}

\proof

Changing the variable, we need for (a) to show the property for the function $x\mapsto e^{\a x} r(e^x)$ on $\rl.$ Observe from the preceding discussion that this function is log-concave and hence unimodal. Therefore, we need to evaluate the limits at $\pm\infty$, which are clearly both zero for $\a\in (0,1)$, zero and one for $\a = 1$ and $\sqrt{\pi/2}$ and zero for $\a = 0.$ This completes the proof of (a). The argument for (b) is the same, starting from the fact that
$$-r'(\sqrt{2x})\; =\; -\sqrt{x}\, \Psi'_x(1/2,1/2,x)\; =\; \sqrt{\frac{x}{\pi}}\,\int_0^\infty e^{-xt} \,\frac{\sqrt{t}}{1+t}\; dt$$  
is also HCM by Property iv) p.68, Theorem 5.4.1, Theorem 4.3.1 and Section 3.5 in \cite{Bond}, since the function $t\mapsto e^{-\varepsilon t}\sqrt{t}/(1+t)$ is up to normalization an HCM density for every $\varepsilon > 0.$  Hence, the function $x\mapsto -e^{\a x} r'(e^x)$ is log-concave and unimodal on $\rl$ for every $\a\in\rl,$ and we can apply the same argument as above.

\endproof

Let us now state, as a consequence of Theorem A, a series of CM properties for Mill's ratio and its derivatives, which were actually the starting point of this paper.

\begin{Propo}
\label{CMMill} 
For every integer $n,$ the function
$$z\;\mapsto\; -\frac{\lpa r^{(n)}(\sqrt{z})\rpa^2}{r^{(2n+1)}(\sqrt{z})}$$
is {\em CM}. For $n= 0,1$ it is also {\em LCM}.
\end{Propo}

\proof

It is clear from (\ref{Herman}) that $H_\nu' = -2\nu H_{-\nu-1}$ and, combining this with the first equality in (\ref{Mill}), we get
$$-\frac{\lpa r^{(n)}(\sqrt{2z})\rpa^2}{r^{(2n+1)}(\sqrt{2z})}\; =\; \frac{(n!)^2\, \lpa H_{-n-1} (\sqrt{z})\rpa^2}{(2n+1)!\, H_{-2n-2}(\sqrt{z})}\cdot$$
The conclusion follows from Corollaries \ref{Para1} and \ref{Para2}.
\endproof

As a by-product of Corollary \ref{Hermit}, let us also mention a new proof of an old and accurate upper bound on Mill's ratio due to Sampford \cite{S}. Other proofs of this classic inequality can be found e.g. in Proposition 3 of \cite{SW} or in Theorem 2.3 of \cite{B}, and we do not claim that our argument is comparatively simple.

\begin{Coro}
\label{Samp}
The mapping $x\mapsto 1/r(x)$ is strictly convex on $\rl.$ In particular, one has
$$r(x)\; < \; \frac{4}{3x +\sqrt{x^2 +8}}$$
for every $x > -1.$
\end{Coro}

\proof

One has
$$\lpa \frac{1}{r}\rpa^{{}_{''}}\! =\; \frac{2(r')^2 - rr''}{r^3}\; >\; 0$$
on $\rl$ since a consequence of the first equality in (\ref{Mill}), the aforementioned fact that $H_{-\nu}' = -2\nu H_{-\nu-1}$, and Remark \ref{DuranDuran}, is
$$\frac{2 (r'(\sqrt{2} x))^2}{r(\sqrt{2} x)\, r''(\sqrt{2} x)}\; =\; \frac{(H_{-2}(x))^2}{H_{-1}(x)\, H_{-3}(x)}\; =\; \frac{(D_{-2}(x))^2}{D_{-1}(x)\, D_{-3}(x)} \; > \; 1$$
for every $x\in\rl.$ To obtain the inequality, observe first from $r' = xr -1$ that 
$$2(r')^2 - rr'' \; = \; (x^2 -1)r^2 - 3xr + 2 \; = \; 2r^2\lpa \frac{1}{r} - \frac{3x + \sqrt{x^2+8}}{4}\rpa \lpa \frac{1}{r} - \frac{3x - \sqrt{x^2+8}}{4}\rpa.$$
The obvious bound $xr<1$ shows next $r(3x - \sqrt{x^2+8}) < 3xr < 4$ on $\rl$ and, putting everything together, we obtain
$$r\; < \; \frac{4}{3x +\sqrt{x^2 +8}}$$
on $(-1,\infty)$ as required.
\endproof

\begin{Rem}
{\em The convexity of $x\mapsto 1/r(x)$ on $\rl^+$ is also a consequence of the case $n=0$ in Proposition \ref{CMMill}.}
\end{Rem}

We next state the following open problem on Mill's ratio and its derivatives, which is a refinement of Proposition \ref{CMMill}. Recall indeed that if a function $f$ is CM, then so is $x\mapsto f(\sqrt{x}).$ The statements are supported by simulations, but we could not find any proof. 

\begin{Conj}
\label{CMMI}
For every integer $n,$ the function
$$z\;\mapsto\; -\frac{\lpa r^{(n)}(z)\rpa^2}{r^{(2n+1)}(z)}$$
is {\em CM}.
\end{Conj}

To handle this last problem, it is natural to start with the representation
$$r^{(p)}(\sqrt{2}\, x)\; =\; (-1)^p\, 2^{\frac{p-1}{2}}\, \Ga((p+1)/2) \,\esp\!\lcr e^{-2x \sqrt{\G_{(p+1)/2}}}\,\rcr$$
for every integer $p$, which follows from (\ref{Herman}) and (\ref{Mill}). In this respect the following proposition, which is interesting in its own right, gives another support for the validity of Conjecture \ref{CMMI}. Unfortunately, we cannot use this factorization directly as we did in Corollaries \ref{B1} and \ref{Para1}, since the factor $\sqrt{\G_{2a}}$ is not ID and hence not self-decomposable.

\begin{Propo}
\label{Half-Gaussian}
For every $a > 0,$ one has 
$$\sqrt{\G_a}\, +\sqrt{\G_a}\; \elaw\; \sqrt{\G_{2a}\times\lpa 1+\sqrt{\B_{a,1/2}}\rpa}.$$
\end{Propo}

\proof
The density of $\sqrt{\G_a}$ is $(2/\Ga(a))\, x^{2a-1} e^{-x^2}$ and its convolution with itself gives the expression
\begin{eqnarray*}
\frac{4x^{4a-1}}{(\Ga(a))^2} \int_0^1 e^{-x^2(y^2+(1-y)^2)} (y(1-y))^{2a-1} dy & = & \frac{8x^{4a-1}}{(\Ga(a))^2} \int_{1/2}^1 e^{-x^2(y^2+(1-y)^2)} (y(1-y))^{2a-1} dy\\
& = & \frac{4^{2-2a} x^{4a-1}}{(\Ga(a))^2} \int_0^1 e^{-x^2(1 + z^2)/2}\, (1-z^2)^{2a-1}\, dz
\end{eqnarray*}
which we can reinject to evaluate the Mellin transform of $\sqrt{\G_a} +\! \sqrt{\G_a}$, by Fubini's theorem: one finds
\begin{eqnarray*}
\esp\lcr \lpa\! \sqrt{\G_a} + \sqrt{\G_a} \rpa^s\rcr & = & \frac{4^{2-2a}}{(\Ga(a))^2} \int_0^\infty\!\!\int_0^1 u^{s+4a-1} e^{-u^2} \lpa\frac{2}{1+z^2}\rpa^{2a+s/2}\!\! (1-z^2)^{2a-1}\, du\, dz\\
& = & \frac{\Ga(2a+s/2)}{\Ga(2a)}\, \times\, \lpa \frac{4^{1-2a}\,\Ga(2a)}{(\Ga(a))^2} \int_0^1 \lpa\frac{2}{1+z^2}\rpa^{2a+s/2} \!\! (1-z^2)^{2a-1}\, dz\rpa\\
& = & \esp\lcr \lpa \sqrt{\G_{2a}}\rpa^{\! s}\rcr\times\, \lpa \frac{2^{1-4a}\,\Ga(2a)}{(\Ga(a))^2} \int_0^1 \lpa\frac{2}{1+z}\rpa^{2a+s/2} \!\! (1-z)^{2a-1}\, z^{-1/2}\, dz\rpa\\
& = & \esp\lcr \lpa \sqrt{\G_{2a}}\rpa^{\! s}\rcr\times\,\esp \lcr \lpa \sqrt{\lpa\frac{2}{1+ \B_{1/2,2a}}\rpa^{(2a)}}\rpa^{\!\! s\,}\rcr
\end{eqnarray*}  
for every $s > -4a.$ Identifying the Mellin factors, we are reduced to show the identity
$$\lpa\frac{2}{1+ \B_{1/2,2a}}\rpa^{(2a)} \elaw\;\, 1\; +\; \sqrt{\B_{a,1/2}}\, ,$$
which follows from the fact that the density of the random variable on both sides is proportional to $(x-1)^{2a-1}(2x-x^2)^{-1/2}$ on $(1,2).$ We omit details.

\endproof

We would like to conclude the paper with the following curious identity, which is a direct consequence of \eqref{base}, Theorem A and Proposition \ref{Half-Gaussian}.

\begin{Coro}
\label{34}
For every $a > 0,$ one has
$$\B'_{a,1/2}\, +\, \B'_{a,1/2}\; \elaw\; \frac{\lpa \sqrt{\G_a} + \sqrt{\G_a}\rpa^2}{\G_{1/2}}\cdot$$
\end{Coro}

\section*{Acknowledgement}
Rui A. C. Ferreira was supported by the ``Funda\c{c}\~{a}o para a Ci\^encia e a Tecnologia (FCT)" through the program ``Stimulus of Scientific Employment, Individual Support-2017 Call" with reference CEECIND/00640/2017.

\end{document}